# Least squares type estimation of the transition density of a particular hidden Markov chain

**Claire Lacour**

*MAP5, Université Paris Descartes*
*45 rue des Saints-Pères*
*75270 Paris Cedex 06*
*France*
*e-mail:* lacour@math-info.univ-paris5.fr

**Abstract:** In this paper, we study the following model of hidden Markov chain: $Y_i = X_i + \varepsilon_i$, $i = 1, \dots, n+1$ with $(X_i)$ a real-valued stationary Markov chain and $(\varepsilon_i)_{1 \leq i \leq n+1}$ a noise having a known distribution and independent of the sequence $(X_i)$. We present an estimator of the transition density obtained by minimization of an original contrast that takes advantage of the regressive aspect of the problem. It is selected among a collection of projection estimators with a model selection method. The $L^2$-risk and its rate of convergence are evaluated for ordinary smooth noise and some simulations illustrate the method. We obtain uniform risk bounds over classes of Besov balls. In addition our estimation procedure requires no prior knowledge of the regularity of the true transition. Finally, our estimator permits to avoid the drawbacks of quotient estimators.

## 1. Introduction

In this paper we consider the following additive hidden Markov model:

$$Y_i = X_i + \varepsilon_i \qquad i = 1, \dots, n+1 \tag{1}$$

with $(X_i)_{i\geq 1}$ a real-valued Markov chain, $(\varepsilon_i)_{i\geq 1}$ a sequence of independent and identically distributed variables and

$$(X_i)_{i\geq 1} \text{ and } (\varepsilon_i)_{i\geq 1} \text{ independent.} \tag{2}$$

Only the variables $Y_1, ..., Y_{n+1}$ are observed. Besides its initial distribution, the chain $(X_i)_{i\geq 1}$ is characterized by its transition, i.e. the distribution of $X_{i+1}$ given $X_i$. We assume that this transition has a density $\Pi$, defined by $\Pi(x, y)dy = P(X_{i+1} \in dy | X_i = x)$, and our aim is to estimate this transition density $\Pi$.

This model belongs to the class of hidden Markov models. The Hidden Markov Models constitute a very famous class of discrete-time stochastic processes, with many applications in various areas such as biology, speech recognition or finance. For a general reference on these models, we refer to Cappé et al. (2005). Here, we study a simple model of HMM where the noise is additive (which allows dealing also with multiplicative noise by use of a logarithm). In standard HMM, it is assumed that the joint density of $(X_i, Y_i)$ has a parametric form and the aim is then to infer the parameter from the observations $Y_1, ..., Y_n$, generally by maximizing the likelihood. For this type of study, see, among others, Baum and Petrie (1966), Leroux (1992), Bakry et al. (1997), Bickel et al. (1998), Jensen and Petersen (1999), Douc et al. (2004), Fuh (2006).

This model is also similar to the so-called convolution model (for which the aim is to estimate the density of $(X_i)_{i\geq 1}$). As in that model, we use the Fourier transform extensively. The restrictions on error distribution and rate of convergence obtained for our estimator are also of the same kind. Related works include Stefanski (1990), Fan (1993), Masry (1993) (for the multivariate case), Pensky and Vidakovic (1999), Comte et al. (2006).

The estimation of the transition density of a hidden Markov chain is studied by Clémençon (2003). His estimator is based on the thresholding of a wavelet-vaguelette decomposition. The drawback is that this estimator does not achieve the minimax rate because of a logarithmic loss. Lacour (2007b) describes an estimation procedure by quotient of an estimator of the joint density and an estimator of the stationary density $f$. The minimax rate is reached by this estimator if it is assumed that $f$ and $f.\Pi$ have the regularity $\alpha$. But this smoothness condition on $f$ raises a problem. Indeed Clémençon (2000) gives an example in which the stationary density $f$ is not continuous, whereas the transition density $\Pi$ is constant. It shows that $f$ can be much less regular than $\Pi$. Therefore, our aim is to find an estimator of the transition density which does not have the above mentioned disadvantages.

To estimate $\Pi$, we use an original contrast inspired by the mean square contrast. The first idea is to connect our problem with the regression model. For any function $g$, we can write

$$g(X_{i+1}) = \left(\int \Pi(.,y)g(y)dy\right)(X_i) + \eta_{i+1}$$

where $\eta_{i+1} = g(X_{i+1}) - \mathbb{E}[g(X_{i+1})|X_i]$. Then, for all function $g$, we can consider $\int \Pi g$ as a regression function. The mean square contrast to estimate this regression function, if the $X_i$ were known, should be $(1/n)\sum_{i=1}^n [t^2(X_i) - 2t(X_i)g(X_{i+1})]$. If $\int g^2 = 1$, this contrast can be written

$$(1/n)\sum_{i=1}^n [\int T^2(X_i, y)dy - 2T(X_i, X_{i+1})]$$

by setting $T(x,y) = t(x)g(y)$ i.e. $T$ such that $\int T(x,y)g(y)dy = t(x)$. It is this contrast which is used in Lacour (2007a) but in our case, only the $Y_1, \ldots, Y_{n+1}$ are known. Therefore we introduce in this paper two operators $Q$ and $V$ such that $\mathbb{E}[Q_{T^2}(Y_i)|X_i] = \int T^2(X_i, y)dy$ and $\mathbb{E}[V_T(Y_i, Y_{i+1})|X_i, X_{i+1}] = T(X_i, X_{i+1})$. It leads to the following contrast:

$$\gamma_n(T) = \frac{1}{n}\sum_{i=1}^{n}[Q_{T^2}(Y_i) - 2V_T(Y_i, Y_{i+1})]. \quad (3)$$

A collection of estimators is then defined by minimization of this contrast on wavelet spaces. Indeed wavelets have many useful properties and in particular they can have a compact support and can be regular enough to balance the smoothness of the noise. A general reference on the subject is Meyer (1990)'s book.

A method of model selection inspired by Barron et al. (1999) and based on contrast (3) is used to build our final estimator. A data driven choice of model is performed via the minimization of a penalized criterion. The chosen model is the one which minimizes the empirical risk added to a penalty function. In most cases when estimating mixing processes, a mixing term appears in this penalty. In the same way, some unknown terms derived from the dependence between the $X_i$ appears in the thresholding constant used to define the estimator of Clémençon (2003). Here a conditioning argument enables to avoid such a mixing term in the penalty. Our penalty contains only known quantities or terms that can be estimated and is then computable.

For an ordinary smooth noise with regularity $\gamma$, the rate of convergence $n^{-\alpha/(2\alpha+4\gamma+2)}$ is obtained if it is assumed that the transition $\Pi$ belongs to a Besov space with regularity $\alpha$. Our estimator is then better than that of Clémençon (2003) which achieves only the rate $(\ln(n)/n)^{\alpha/(2\alpha+4\gamma+2)}$. Moreover this rate is obtained without assuming the regularity $\alpha$ of $\Pi$ known.

This paper is organized as follows. In Section 2 we present the model and the assumptions. Section 3 is devoted to the definitions of the contrast and of the estimator. The main result and a sketch of proof are to be found in Section 4. Numerical illustrations through simulated examples are reported in Section 5. The detailed proofs are gathered in Section 6.

## 2. Study framework

### 2.1. Notations

For the sake of clarity, we use lowercase letters for dimension 1 and capital letters for dimension 2. For a function $t : \mathbb{R} \mapsto \mathbb{R}$, we denote by $\|t\|$ the $L^2$ norm that is $\|t\|^2 = \int_{\mathbb{R}} t^2(x)dx$. The Fourier transform $t^*$ of $t$ is defined by

$$t^*(u) = \int e^{-ixu}t(x)dx.$$

Notice that the function $t$ is the inverse Fourier transform of $t^*$ and can be written $t(x) = 1/(2\pi)\int e^{ixu}t^*(u)du$. The convolution product is defined by $(t * s)(x) = \int t(x-y)s(y)dy$.

In the same way, for a function $T : \mathbb{R}^2 \mapsto \mathbb{R}$, $\|T\|^2 = \iint_{\mathbb{R}^2} T^2(x,y)dxdy$ and

$$T^*(u,v) = \iint e^{-ixu-iyv}T(x,y)dxdy.$$

We denote by $t \otimes s$ the function: $(x,y) \mapsto (t \otimes s)(x,y) = t(x)s(y)$.

We will estimate $\Pi$ on a compact set $A = A_1 \times A_2$ only and we denote by $\|.\|_A$ the norm in $L^2(A)$ i.e.

$$\|T\|_A^2 = \iint_A T^2(x,y)dxdy.$$

### *2.2. Assumptions on noise*

The Markov chain $(X_i)_{i\geq 1}$ is observed through a noise sequence $(\varepsilon_i)_{i\geq 1}$ of independent and identically distributed random variables. The density of $\varepsilon_i$ is denoted by $q$ and is assumed to be known. We assume that the Fourier transform of $q$ never vanishes and that $q$ is ordinary smooth. More precisely the assumption on the error density is the following:

**H1** $q$ is uniformly bounded and there exist $\gamma > 0$ and $k_0 > 0$ such that $\forall x \in \mathbb{R}$ $|q^*(x)| \geq k_0(x^2+1)^{-\gamma/2}$.

This assumption restrains the regularity class of the noise. Among the so-called ordinary smooth noises, we can cite the Laplace distribution, the exponential distribution and all the Gamma or symmetric Gamma distributions. The noise follows a Gamma distribution with scale parameter $\lambda$ and shape parameter $\zeta$ if $q(x) = \lambda^\zeta x^{\zeta-1}e^{-\lambda x}/\Gamma(\zeta)$ for $x > 0$ with $\Gamma$ the classic Gamma function. Then

$$|q^*(x)| = \left(1 + \frac{x^2}{\lambda^2}\right)^{-\zeta/2}.$$

So $q$ is bounded and verifies H1 with $\gamma = \zeta$. The case $\zeta = 1$ corresponds to an exponential distribution and if $\lambda = 1/2$, $\zeta = p/2$, it is a chi-square $\chi(p)$. A Laplace noise is defined in the following way

$$q(x) = \frac{\lambda}{2}e^{-\lambda|x-\mu|} \text{ and } |q^*(x)| = \frac{\lambda^2}{x^2+\lambda^2}$$

Then H1 is satisfied with $\gamma = 2$. More generally, we can define the symmetric gamma distribution with density $q(x) = \lambda^\zeta|x|^{\zeta-1}e^{-\lambda|x|}/(2\Gamma(\zeta))$. The characteristic function is then

$$q^*(x) = \left(1 + \frac{x^2}{\lambda^2}\right)^{-\zeta/2} \cos\left(2\zeta \arctan\left(\frac{x}{\lambda + \sqrt{x^2+\lambda^2}}\right)\right)$$

so that H1 is verified with $\gamma = \zeta + 1$ if $\zeta$ is an odd integer and $\gamma = \zeta$ otherwise.

**Remark 1.** *We have to point out that the Gaussian noise does not verify Assumption H1. Indeed, an exponential decrease of the Fourier transform of the error density is more difficult to control and a supersmooth noise makes denoising more difficult. For that reason, many authors, among which Butucea (2004), Koo and Lee (1998) or Youndjé and Wells (2002), have considered only ordinary smooth noise. The method used in this paper does not allow dealing with supersmooth noise. Indeed, it requires a wavelet basis more regular than the noise and with compact support (because of Assumption H4 below), which is impossible when the noise is supersmooth.*

### *2.3. Assumptions on the chain*

The hypotheses on the hidden Markov chain $(X_i)_{i\geq 1}$ are the following:

**H2** The chain is irreducible, positive recurrent and stationary with unknown density $f$.

**H3** There exists a positive real $f_0$ such that, for all $x$ in $A_1$,

$$f_0 \leq f(x) \leq \|f\|_{\infty,A_1} < \infty$$

**H4** The transition density $\Pi$ is bounded on $A$ by $\|\Pi\|_{\infty,A} < \infty$.

**H5** The process $(X_k)$ is geometrically $\beta$-mixing ($\beta_q \leq e^{-\theta q}$), or arithmetically $\beta$-mixing ($\beta_q \leq q^{-\theta}$) with $\theta > 8$ where

$$\beta_q = \int \|P^q(x,.) - \mu\|_{TV} f(x)dx$$

with $P^q(x,.)$ the distribution of $X_{i+q}$ given $X_i = x$, $\mu$ the stationary distribution and $\|.\|_{TV}$ the total variation distance.

We refer to Doukhan (1994) for details on the $\beta$-mixing. Assumption H5 implies that the process $(Y_k)$ is $\beta$-mixing, with $\beta$-mixing coefficients smaller than those of $(X_k)$. Assumption H3 is common (but restrictive) and is crucial to control the empirical processes brought into play. A lot of processes verify Assumptions H2–H5, as autoregressive processes, diffusions or ARCH processes. These examples are detailed in Lacour (2007a).

## 3. Estimation procedure

### *3.1. Projection spaces*

Here we describe the projection that we use to estimate the transition $\Pi$. We will consider an increasing sequence of spaces, indexed by $m$, to construct a collection of estimators. For the sake of simplicity, we assume that $A = [0,1]^2$.

We use a compactly supported wavelet basis on the interval $[0,1]$, described in Cohen et al. (1993). The construction provides a set of functions $(\phi_k)$ for $k = 0,\dots,2^J - 1$ with $J$ a fixed level, and for all $j > J$ a set of functions $(\psi_{jk}), k = 0, \dots, 2^j - 1$. The collection of these functions forms a complete orthonormal system on $[0,1]$. Then , for $u$ in $L^2([0,1])$, we can write

$$u = \sum_{k=0}^{2^J-1} b_k \phi_k + \sum_{j>J} \sum_{k=0}^{2^j-1} a_{jk} \psi_{jk}.$$

Actually

$$\phi_k(x) = \begin{cases} 2^{J/2}\phi^0(2^J x - k) & \text{if } k = 0, \dots, N-1 \\ 2^{J/2}\phi(2^J x - k) & \text{if } k = N, \dots, 2^J - N - 1 \\ 2^{J/2}\phi^1(2^J x - k) & \text{if } k = 2^J - N, \dots, 2^J - 1 \end{cases}$$

where $\phi$ is a Daubechies father wavelet with support $[-N+1, N]$ and $\phi^0$, $\phi^1$ are edge wavelets explicitly constructed in Cohen et al. (1993). The functions $\phi_k$ have support $[(k-N+1)/2^J, (k+N)/2^J] \cap [0,1]$. For $r$ a positive real, $N$ is chosen large enough so that $\phi$ has regularity $r$ (in the sense defined in (4)): this is possible since it is a property of the Daubechies wavelets that the smoothness of $\phi$ increases linearly with $N$. We choose $J$ such that $2^J \geq 2N$ so that the two edges do not interact (no overlap between $\phi^0$ and $\phi^1$). The construction ensures that $\phi^0$ and $\phi^1$ are also of regularity $r$. In the same way, for each level $j$, the $\psi_{jk}$ are dilatation and translation of functions $\psi$, $\psi^0$ and $\psi^1$ with regularity $r$.

Now we construct a wavelet basis of $L^2([0,1]^2)$ by the tensorial product method (see Meyer (1990) Chapter 3 Section 3). The father wavelet is $\phi \otimes \phi$ and the mother wavelets are $\phi \otimes \psi$, $\psi \otimes \phi$, $\psi \otimes \psi$. A function $T$ in $L^2([0,1]^2)$ can then be written

$$T(x,y) = \sum_{k=0}^{2^J-1} \sum_{l=0}^{2^J-1} b_{kl} \phi_k(x) \phi_l(y) + \sum_{j>J} \sum_{k=0}^{2^j-1} \sum_{l=0}^{2^j-1} (a_{jkl}^{(1)} \phi_{jk}(x) \psi_{jl}(y) + a_{jkl}^{(2)} \psi_{jk}(x) \phi_{jl}(y) + a_{jkl}^{(3)} \psi_{jk}(x) \psi_{jl}(y)).$$

For the sake of simplicity, we adopt the following notation

$$T(x,y) = \sum_{j \geq J} \sum_{(k,l) \in \Lambda_j} a_{jkl} \varphi_{jk}(x) \varphi_{jl}(y).$$

where $\varphi_{jk} = 2^{j/2} \varphi(2^j x - k)$ with $\varphi = \phi, \phi^0, \phi^1, \psi, \psi^0$ or $\psi^1$ according to the values of $j$ and $k$. For $j > J$, $\Lambda_j$ is a set with cardinal $3.2^{2j}$ and $\Lambda_J$ is a set with cardinal $2^{2J}$. In the rest of this paper we will use the following property of $\varphi$ deriving from the regularity of the initial Daubechies wavelet: there exists a positive constant $k_3$ such that

$$\forall u \in \mathbb{R} \quad |\varphi^*(u)| \leq k_3 (u^2+1)^{-r/2} \tag{4}$$

Now , for $m \geq J$, we can consider the space

$$\mathbb{S}_m = \{T : \mathbb{R}^2 \to \mathbb{R}, \quad T(x,y) = \sum_{j=J}^{m} \sum_{(k,l) \in \Lambda_j} a_{jkl} \varphi_{jk}(x) \varphi_{jl}(y)\}$$

Note that the functions in $\mathbb{S}_m$ are all supported in the interval $[0,1]^2$. The dimension of the space $\mathbb{S}_m$ is $D_m^2 = 2^{2J} + 3\sum_{j=J+1}^{m} 2^{2j} \in [2^{2m}, 2^{2m+2}]$. We denote by $\mathcal{S}$ the space $S_{m_0}$ with the greatest dimension $D_{m_0}^2 = \mathcal{D}^2$ smaller than $n^{1/(4\gamma+2)}$. It is the maximal space that we consider. The spaces $\mathbb{S}_m$ have the following properties:

**P1** $m' \leq m \Rightarrow \mathbb{S}_{m'} \subset \mathbb{S}_m$
**P2** $\| \sum_{jkl} a_{jkl}\varphi_{jk} \otimes \varphi_{jl}\|^2 = \sum_{jkl} a_{jkl}^2$.

This property derives from the orthonormality of the basis.

Now, for all function $t : \mathbb{R} \mapsto \mathbb{R}$, let $v_t$ be the inverse Fourier transform of $t^*/q^*(-.)$, i.e.

$$v_t(x) = \frac{1}{2\pi}\int e^{ixu}\frac{t^*(u)}{q^*(-u)}du.$$

This operator is introduced because it verifies $\mathbb{E}[v_t(Y_k)|X_k] = t(X_k)$ for all function $t$. We can write the following lemma :

**Lemma 1.** *If $r > \gamma + 2$, there exists $\Phi_1 > 0$ such that*

**P3** $\| \sum_{j=J}^{m}\sum_k \varphi_{jk}^2\|_\infty \leq \Phi_1 D_m$
**P4** $\| \sum_k |v_{\varphi_{jk}}|^2\|_\infty \leq \Phi_1(2^j)^{2\gamma+2}$
**P5** $\sum_k \|v_{\varphi_{jk}}\|^2 \leq \Phi_1(2^j)^{2\gamma+1}$
**P6** $\| \sum_{kk'} |v_{\varphi_{jk}\varphi_{jk'}}|^2\|_\infty \leq \Phi_1(2^j)^{2\gamma+3}$
**P7** $\sum_{kk'} \int |v_{\varphi_{jk}\varphi_{jk'}}|^2 \leq \Phi_1(2^j)^{2\gamma+2}$

This lemma is proved in Section 6.

### *3.2. Construction of a contrast*

Now let us estimate the transition density of the Markov chain by minimizing a contrast. This section is devoted to the definition of this contrast. We explain here how it can be obtained, first by considering the case without noise.

#### *3.2.1. First step: if $X_1, \dots, X_{n+1}$ were observed*

We present here a heuristic to understand why we choose the contrast, assuming first that the $(X_i)$ are known. For all function $g$, the definition of the transition density implies $\mathbb{E}[g(X_{i+1})|X_i] = \int \Pi(X_i, y)\; g(y)dy$ so that we can write

$$g(X_{i+1}) = \left(\int \Pi(., y)g(y)dy\right)(X_i) + \eta_i$$

where $\eta_i = g(X_{i+1}) - \mathbb{E}[g(X_{i+1})|X_i]$ is a centered process. We recognize then a regression model. A contrast to estimate $\int \Pi(., y)g(y)dy$ is

$$\gamma_n(u) = \frac{1}{n}\sum_{i=1}^{n}[u^2(X_i) - 2u(X_i)g(X_{i+1})].$$

It is the classical mean square contrast to estimate a regression function. But we want to estimate $\Pi(.,y)$ and not only $\int \Pi(.,y)g(y)dy$.

So we observe that if $\int g^2=1$ and $T(x,y)=u(x)g(y)$, then $u(.)=\int T(.,y)g(y)dy$. So if $u(.)=\int T(.,y)g(y)dy$ estimates $\int \Pi(.,y)g(y)dy$, we can assume that $T$ estimates $\Pi$. Since $\int T^2(.,y)dy = u^2(.)$, the contrast becomes

$$\gamma_n(T) = \frac{1}{n}\sum_{i=1}^{n}[\int T^2(X_i,y)dy - 2T(X_i,X_{i+1})]$$

It is the contrast studied in Lacour (2007a) and it allows for a good estimation of $\Pi(.,y)$ when the Markov chain is observed. We can observe that

$$\mathbb{E}\gamma_n(T) = \int T^2(x,y)f(x)dxdy - 2\int T(x,y)f(x)\Pi(x,y)dxdy = \|T-\Pi\|_f^2 - \|\Pi\|_f^2$$

where $f$ is the density of $X_i$ and

$$\|T\|_f = \left(\int T^2(x,y)f(x)dxdy\right)^{1/2}.$$

Then this contrast is an empirical counterpart of the distance $\|T-\Pi\|_f$.

### *3.2.2. Second step: the $X_i$'s are unknown, the observations are the $Y_i$'s*

The aim of this step is to modify the previous contrast, to take into account that the $X_i$'s are not observed. To do this, we use the same technique as in the convolution problem (see Comte et al. (2006)). Let us denote by $F_X$ the density of $(X_i,X_{i+1})$ and $F_Y$ the density of $(Y_i,Y_{i+1})$. We remark that $F_Y = F_X*(q\otimes q)$ and $F_Y^* = F_X^*(q^*\otimes q^*)$ and then

$$\mathbb{E}[T(X_i,X_{i+1})] = \iint TF_X = \frac{1}{2\pi}\iint T^*\overline{F_X^*} = \frac{1}{2\pi}\iint \frac{T^*}{\overline{q^*}\otimes\overline{q^*}}\overline{F_Y^*}$$

by using the Parseval equality. The idea is then to define $V_T^* = T^*/(\overline{q^*}\otimes\overline{q^*})$ so that

$$\mathbb{E}[T(X_i,X_{i+1})] = \frac{1}{2\pi}\iint V_T^*\overline{F_Y^*} = \iint V_T F_Y = \mathbb{E}[V_T(Y_i,Y_{i+1})].$$

Then we replace the term $T(X_i,X_{i+1})$ in the contrast by $V_T(Y_i,Y_{i+1})$. In the same way, we find an operator $Q$ to replace the term $\int T^2(X_i,y)dy$. More precisely, for all function $T$, let $V_T$ be the inverse Fourier transform of $T^*/(q^*\otimes q^*)(-.)$, i.e.

$$V_T(x,y) = \frac{1}{4\pi^2}\iint e^{ixu+iyv}\frac{T^*(u,v)}{q^*(-u)q^*(-v)}dudv.$$

Let $Q_T$ be the inverse Fourier transform of $T^*(.,0)/(q^*)(-.)$, i.e.

$$Q_T(x) = \frac{1}{2\pi}\int e^{ixu}\frac{T^*(u,0)}{q^*(-u)}du.$$

$V$ and $Q$ have been chosen so that the following lemma holds.

**Lemma 2.** *For all* $k \in \{1, \dots, n+1\}$

1. $\mathbb{E}[V_T(Y_k, Y_{k+1})|X_1, ..., X_{n+1}] = T(X_k, X_{k+1})$
2. $\mathbb{E}[V_T(Y_k, Y_{k+1})] = \iint T(x,y)\Pi(x,y)f(x)dxdy$
3. $\mathbb{E}[Q_T(Y_k)|X_1, ..., X_{n+1}] = \int T(X_k, y)dy$
4. $\mathbb{E}[Q_T(Y_k)] = \iint T(x,y)f(x)dxdy$

Points *1.* and *3.* are proved in Section 6, the other assertions are their immediate consequences. Note that $V$ and $Q$ are strongly linked with $v$. In particular $V_{s\otimes t}(x,y) = v_s(x)v_t(y)$ and $Q_{s\otimes t}(x) = v_s(x)\int t(y)dy$.

By using the operators $V$ and $Q$, we now define the contrast, depending only on the observations $Y_1, \dots, Y_{n+1}$:

$$\gamma_n(T) = \frac{1}{n}\sum_{k=1}^{n}[Q_{T^2}(Y_k) - 2V_T(Y_k, Y_{k+1})]$$

With Lemma 2, we compute $\mathbb{E}(\gamma_n(T)) = \iint T^2(x,y)f(x)dxdy - 2\iint T(x,y)$ $\Pi(x,y)f(x)dxdy = \|T-\Pi\|_f^2 - \|\Pi\|_f^2$. So we want to estimate $\Pi$ by minimizing $\gamma_n$. The definition of the contrast leads to the following "empirical norm":

$$\Psi_n(T) = \frac{1}{n}\sum_{k=1}^{n} Q_{T^2}(Y_k).$$

The term empirical norm is used because $\mathbb{E}\Psi_n(T) = \|T\|_f^2$, but $\Psi_n$ is not a norm in the common sense of the word.

### *3.3. Definition of the estimator*

We have to minimize the contrast $\gamma_n$ to find our estimator. By writing $T = \sum_{j=J}^{m}\sum_{(k,l)\in\Lambda_j} a_{jkl}\varphi_{jk}\otimes\varphi_{jl} = \sum_\lambda a_\lambda\ \omega_\lambda(x,y)$, we obtain

$$\frac{\partial\gamma_n(T)}{\partial a_{\lambda_0}} = \frac{2}{n}\sum_{i=1}^{n}\left(\sum_\lambda a_\lambda Q_{\omega_\lambda\omega_{\lambda_0}}(Y_i) - V_{\omega_{\lambda_0}}(Y_i, Y_{i+1})\right).$$

Then, by denoting $A_m$ the vector of the coefficients $a_\lambda$ of $T$,

$$\forall\lambda_0 \quad \frac{\partial\gamma_n(T)}{\partial a_{\lambda_0}} = 0 \quad \Longleftrightarrow \quad G_m A_m = Z_m \tag{5}$$

where

$$G_m = \left[\frac{1}{n}\sum_{i=1}^{n} Q_{\omega_\lambda\omega_\mu}(Y_i)\right]_{\lambda,\mu}, \quad Z_m = \left[\frac{1}{n}\sum_{i=1}^{n} V_{\omega_\lambda}(Y_i, Y_{i+1})\right]_\lambda$$

But the matrix $G_m$ is not necessarily invertible. This is why we introduce the set

$$\Gamma = \left\{\min \mathrm{Sp}(G_m) \geq \frac{2}{3}f_0\right\} \tag{6}$$

where Sp denotes the spectrum, i.e. the set of the eigenvalues of the matrix and $f_0$ is the lower bound of $f$ on $A_1$. On $\Gamma$, $G_m$ is invertible and $\gamma_n$ is convex so that the minimization of $\gamma_n$ is equivalent to Equation (5) and admits the solution $A_m = G_m^{-1} Z_m$. Now we can define

$$\hat{\Pi}_m = \begin{cases} \arg\min_{T\in\mathbb{S}_m} \gamma_n(T) & \text{on } \Gamma, \\ 0 & \text{on } \Gamma^c. \end{cases}$$

**Remark 2.** *The term* $2/3$ *in* $\Gamma$ *can be replaced by any constant smaller than 1. Moreover, the construction of* $\hat{\Pi}_m$ *described here requires the knowledge of* $f_0$*. Nevertheless, when* $f_0$ *is unknown, we can replace it by an estimator* $\hat{f}_0$ *defined as the minimum of an estimator of* $f$ *(for an estimator of the density of a hidden Markov chain, see Lacour (2007b)). The result is then unchanged if* $f$ *is regular enough and the mixing rate high enough.*

We then have an estimator of $\Pi$ for all $\mathbb{S}_m$. But we have to choose the best model $m$ to obtain an estimator which achieves the best rate of convergence, whatever the regularity of $\Pi$. So we set

$$\hat{m} = \arg\min_{m\in\mathcal{M}_n} \{\gamma_n(\hat{\Pi}_m) + \text{pen}(m)\}$$

where pen is a penalty function to be specified later and

$$\mathcal{M}_n = \{m \geq J, D_m^{4\gamma+2} \leq n\}.$$

Then we can define our final estimator:

$$\tilde{\Pi} = \begin{cases} \hat{\Pi}_{\hat{m}} & \text{if } \|\hat{\Pi}_{\hat{m}}\| \leq k_n \quad \text{with } k_n = n^{1/2}, \\ 0 & \text{else.} \end{cases}$$

## 4. Result

### 4.1. Risk and rate of convergence

For a function $G$ and a subspace $\mathbb{S}$, we define

$$d_A(G, \mathbb{S}) = \inf_{T\in\mathbb{S}} \|G - T\|_A.$$

We recall that $A$ is the estimation area. For each estimator $\hat{\Pi}_m$, we have the following decomposition of the risk:

**Proposition 1.** *We consider a Markov chain and a noise satisfying Assumptions H1–H5 with* $\gamma \geq 3/4$*. For* $m$ *fixed in* $\mathcal{M}_n$*, we consider* $\hat{\Pi}_m$ *the estimator of the transition density* $\Pi$*, previously described. Then there exists* $C > 0$ *such that*

$$\mathbb{E}\|\hat{\Pi}_m - \Pi\|_A^2 \leq C \left\{ d_A^2(\Pi, S_m) + \frac{D_m^{4\gamma+2}}{n} \right\}.$$

We do not prove this proposition because this result is included in Theorem 1 below, which is proved in Section 6.

Now if $\Pi$ belongs to a Besov space with regularity $\alpha$, it is a common approximation property of the wavelet spaces that $d_A^2(\Pi, S_m) \leq C D_m^{-2\alpha}$. So, choosing $m_1$ such that $D_{m_1} = n^{1/(2\alpha+4\gamma+2)}$, we obtain the minimum risk

$$\mathbb{E}\|\hat{\Pi}_{m_1} - \Pi\|_A^2 \leq C n^{-\frac{2\alpha}{2\alpha+4\gamma+2}}.$$

But this choice of $m_1$ is impossible if $\alpha$ is unknown (it is *a priori* the case since $\Pi$ is unknown). That is why we have built our estimator $\tilde{\Pi}$ via model selection. Now we can state the following theorem.

**Theorem 1.** *We consider a Markov chain and a noise that satisfy Assumptions H1–H5 with $\gamma > 3/4$. We consider $\tilde{\Pi}$ the estimator of the transition density $\Pi$ previously described with $r > 2\gamma + 3/2$ and*

$$\text{pen}(m) = K \frac{D_m^{4\gamma+2}}{n} \quad \textit{for some } K > K_0$$

*where $K_0 = C(\gamma)\Phi_1^2 \|q\|_\infty^2 f_0^{-1}$. Then there exists $C' > 0$ such that*

$$\mathbb{E}\|\tilde{\Pi} - \Pi\|_A^2 \leq C \inf_{m \in \mathcal{M}_n} \{d_A^2(\Pi, S_m) + \text{pen}(m)\} + \frac{C'}{n}$$

*with $C = \max(2 + 72 f_0^{-1} \|f\|_{\infty, A_1}(1 + 2\|\Pi\|_A^2), 12 f_0^{-1}(1 + 2\|\Pi\|_A^2))$.*

Note that this result is non-asympotic. It is an advantage of the least square method over the quotient method.

All the constants on which the penalty depends do not have the same status. The constants $\Phi_1$, $\gamma$ and $\|q\|_\infty$ are known, since the wavelet basis and the noise distribution are known. The constant $f_0$ is unknown but it can be estimated (see Remark 2). Then, even if it means replacing $f_0$ by an estimator $\hat{f}_0$, the penalty is computable. In particular the dependence coefficients of the sequence do not appear at all in the penalty.

The condition $\gamma > 3/4$ is due to an additional term of order $D_m^{2\gamma+7/2}/n$ (coming from the term $(1/n)\sum_{i=1}^n Q_{T^2}(Y_i)$ in the contrast) inside the penalty. If $\gamma > 3/4$, then $2\gamma + 7/2 < 4\gamma + 2$ and $D_m^{4\gamma+2}/n$ is the dominant term. If $\gamma = 3/4$, the result is still true but the constant in the penalty also depends on $\|\Pi\|_A$. In the other cases the estimation is possible but the term $D_m^{2\gamma+7/2}/n$ is not negligible any more and the order of the variance (and consequently the rate of convergence) must be changed. This constraint $\gamma > 3/4$ is not very restrictive since $\gamma$ must be larger than $1/2$ in order that $q$ be square integrable. Moreover in the case of a Gamma noise, $q$ is not bounded if $\gamma < 1$.

We can now evaluate the rate of convergence of our estimator.

**Corollary 1.** *We suppose that the restriction of $\Pi$ to $A$ belongs to the Besov space $B_{2,\infty}^\alpha(A)$ with $\alpha < r$. Then, under the assumptions of Theorem 1,*

$$\mathbb{E}\|\tilde{\Pi} - \Pi\|_A^2 = O(n^{-\frac{2\alpha}{2\alpha+4\gamma+2}}).$$

To our knowledge, the minimax rates are unknown in the specific estimation problem we consider here and finding them is definitely beyond the scope of this paper. Nevertheless, Clémençon (2003) proved that the rate $n^{-\frac{2\alpha}{2\alpha+4\gamma+2}}$ is optimal whenever $f$ and $f\Pi$ belong to $B^{\alpha}_{2,\infty}(\mathbb{R})$ and $B^{\alpha}_{2,\infty}(\mathbb{R}^2)$ respectively.

Nevertheless we remark that we obtain the same rate of convergence with $\tilde{\Pi}$ as those obtained with $\hat{\Pi}_{m_1}$ where $D_{m_1} = n^{1/(2+4\gamma+2\alpha)}$, but without requiring the knowledge of $\alpha$. Moreover our estimator is better than the one of Clémençon (2003), which achieves only the rate $(\ln(n)/n)^{\frac{2\alpha}{2\alpha+4\gamma+2}}$. It is also an improvement on the result of Lacour (2007b) because this rate is obtained without requiring any regularity for $f$ or $f\Pi$.

If we want to compare the quotient method described in Lacour (2007b) and the one introduced in this paper, we can say that only the quotient method allows dealing with supersmooth distributions, at least from a theoretical point of view. However, the least squares method has the advantage of giving a good rate of convergence without requiring prior information on the stationary density. Moreover, our result is non-asymptotic contrary to the one of Lacour (2007b).

### 4.2. Sketch of proof of Theorem 1

We give in this section a sketch of proof of Theorem 1.

Let $m \in \mathcal{M}_n$. We denote by $\Pi_m$ the orthogonal projection of $\Pi$ on $\mathbb{S}_m$. We have the following bias-variance decomposition

$$\mathbb{E}\|\tilde{\Pi} - \Pi\|_A^2 = \mathbb{E}\|\tilde{\Pi} - \Pi_m\|_A^2 + \|\Pi_m - \Pi\|_A^2$$

The term $\|\tilde{\Pi} - \Pi_m\|_A^2$ can be written in the following way

$$\begin{aligned}\|\tilde{\Pi} - \Pi_m\|_A^2 &= \|\tilde{\Pi} - \Pi_m\|_A^2 \mathbb{1}_{\{\|\hat{\Pi}_{\hat{m}}\| \le k_n\}} + \|\tilde{\Pi} - \Pi_m\|_A^2 \mathbb{1}_{\{\|\hat{\Pi}_{\hat{m}}\| > k_n\}} \\ &\le \|\hat{\Pi}_{\hat{m}} - \Pi_m\|_A^2 + \|\Pi_m\|_A^2 \mathbb{1}_{\{\|\hat{\Pi}_{\hat{m}}\| > k_n\}}\end{aligned}$$

since $\tilde{\Pi} = 0$ on the set $\{\|\hat{\Pi}_{\hat{m}}\| > k_n\}$ and $\tilde{\Pi} = \hat{\Pi}_{\hat{m}}$ on the complement. The term $\|\Pi_m\|_A^2 \mathbb{1}_{\{\|\hat{\Pi}_{\hat{m}}\| > k_n\}}$ is easily dealt with, the main term is $\|\hat{\Pi}_{\hat{m}} - \Pi_m\|_A^2$. But, on $\Gamma$, the definitions of $\hat{\Pi}_m$ and $\hat{m}$ lead to the inequality

$$\gamma_n(\hat{\Pi}_{\hat{m}}) + \mathrm{pen}(\hat{m}) \le \gamma_n(\Pi_m) + \mathrm{pen}(m). \tag{7}$$

Letting $Z_{n,m}(T) = \frac{1}{n}\sum_{k=1}^{n}[V_T(Y_k, Y_{k+1}) - Q_{T\Pi_m}(Y_k)]$, a fast computation gives

$$\gamma_n(\hat{\Pi}_{\hat{m}}) - \gamma_n(\Pi_m) = \Psi_n(\hat{\Pi}_{\hat{m}} - \Pi_m) - 2Z_{n,m}(\hat{\Pi}_{\hat{m}} - \Pi_m)$$

so that (7) becomes

$$\begin{aligned}\Psi_n(\hat{\Pi}_{\hat{m}} - \Pi_m) &\le 2Z_{n,m}(\hat{\Pi}_{\hat{m}} - \Pi_m) + \mathrm{pen}(m) - \mathrm{pen}(\hat{m}) \\ &\le 2\|\hat{\Pi}_{\hat{m}} - \Pi_m\|_f \sup_{T \in B_f(m,\hat{m})} Z_{n,m}(T) + \mathrm{pen}(m) - \mathrm{pen}(\hat{m})\end{aligned}$$

where $B_f(m, \hat{m}) = \{T \in S_m + S_{\hat{m}}, \quad \|T\|_f = 1\}$. The main steps of the proof are then

1. to control the term $\sup_{T\in B_f(m,\hat{m})} Z_{n,m}(T)$,
2. to link the empirical "norm" $\Psi_n$ with the $L^2$ norm $\|.\|_A$.

• To deal with the supremum of the empirical process $Z_{n,m}(T)$, we use an inequality of Talagrand stated in Lemma 6 (Section 6.8). This inequality is very powerful but can be applied only to a sum of independent random variables. That is why we split $Z_{n,m}(T)$ into three processes plus a bias term.

$$Z_{n,m}(T) = Z_n^{(1)}(T) - Z_n^{(2)}(T) + Z_n^{(3)}(T) + \iint T(x,y)(\Pi - \Pi_m)(x,y)f(x)dxdy$$

with

$$\begin{cases} Z_n^{(1)}(T) = \dfrac{1}{n}\displaystyle\sum_{k=1}^{n} V_T(Y_k, Y_{k+1}) - \mathbb{E}[V_T(Y_k, Y_{k+1})|X_1,\dots,X_{n+1}] \\ Z_n^{(2)}(T) = \dfrac{1}{n}\displaystyle\sum_{k=1}^{n} Q_{T\Pi_m}(Y_k) - \mathbb{E}[Q_{T\Pi_m}(Y_k)] \\ Z_n^{(3)}(T) = \dfrac{1}{n}\displaystyle\sum_{k=1}^{n} T(X_k, X_{k+1}) - \mathbb{E}[T(X_k, X_{k+1})] \end{cases}$$

For the first process $Z_n^{(1)}$, we are back to independent variables by remarking that, conditionally to $X_1,\dots,X_{n+1}$, the couples $(Y_{2i-1}, Y_{2i})$ are independent (see Proposition 3).

For the other processes, we use the mixing assumption H5 to build auxiliary variables $X_i^*$ which are approximations of the $X_i$'s and which constitute independent clusters of variables (see Proposition 4).

• To pass from $\Psi_n$ to the $L^2$ norm, we introduce the following set

$$\Delta = \{\forall T \in \mathcal{S} \quad \|T\|_f^2 \le \frac{3}{2}\Psi_n(T)\}$$

We can easily prove (see Section 6.3) that $\Delta \subset \Gamma$. Then,

$$\|\hat{\Pi}_{\hat{m}} - \Pi_m\|_A \mathbb{1}_\Delta \le \frac{3}{2} f_0^{-1} \Psi_n(\hat{\Pi}_{\hat{m}} - \Pi_m)\mathbb{1}_\Gamma$$

It remains to prove that $P(\Delta^c) = P(\exists T \in \mathcal{S}, \Psi_n(T) < (2/3)\mathbb{E}[\Psi_n(T)])$ is small enough. It is done in Proposition 2.

## 5. Simulations

To illustrate the method, we compute our estimator $\tilde{\Pi}$ for different Markov processes with known transition density. The estimation procedure contains several Fourier transforms. This may seem heavy, but, for each noise distribution, the computation of $v_{\varphi_{jk}}$ for all the basis functions can be done beforehand. Here we use the Daubechies wavelet D20. Next, to compute $\tilde{\Pi}$ from data $Y_1,\dots,Y_{n+1}$, we use the following steps (see Section 3.3):

- For each $m$, compute matrices $G_m$ and $Z_m$,
- Deduce the matrix $A_m$,
- Select the $\hat{m}$ which minimizes $\gamma_n(\hat{\Pi}_m) + \text{pen}(m) = -{}^t A_m Z_m + \text{pen}(m)$,
- Compute $\tilde{\Pi}$ using matrix $A_{\hat{m}}$.

Actually, following the theoretical procedure, we should set $\hat{\Pi}_m = 0$ on $\Gamma^c$ (see Section 3.3) but, for practical purposes, it is more sensible to inverse $G_m$ whenever possible. In all the examples examined below, the minimum of the spectrum of $G_m$ has never been too small (so that we merely inverted it without using set $\Gamma$). The reason is that $P(\Gamma^c)$ is very small: it appears in the proofs that it can be bounded with an exponential inequality.

We consider several kinds of Markov chains :

- An autoregressive process denoted by AR and defined by:

$$X_{n+1} = aX_n + b + \varepsilon_{n+1}$$

  where the $\varepsilon_{n+1}$ are independent and identical distributed random variables, with centered Gaussian distribution with variance $\sigma^2$. For this process, the transition density can be written $1/(\sigma\sqrt{2\pi})\exp(-(y - ax - b)^2/2\sigma^2)$. We consider the following parameter values :

  (i) $a = 2/3$, $b = 0$, $\sigma^2 = 5/9$, estimated on $[-2, 2]^2$.

  (ii) $a = 0.5$, $b = 3$, $\sigma^2 = 1$, and then the process is estimated on $[4, 8]^2$.

- A radial Ornstein-Uhlenbeck process (in its discrete version). For $j = 1, \ldots, \delta$, we define the processes: $\xi^j_{n+1} = a\xi^j_n + \beta\varepsilon^j_n$ where the $\varepsilon^j_n$ are i.i.d. standard Gaussian. The chain is then defined by $X_n = \sqrt{\sum_{i=1}^{\delta}(\xi^i_n)^2}$. The transition density is given in Chaleyat-Maurel and Genon-Catalot (2006) where this process is studied in detail:

$$\pi(x, y) = \mathbb{1}_{y>0} \exp\left(-\frac{y^2 + a^2x^2}{2\beta^2}\right) I_{\delta/2-1}\left(\frac{axy}{\beta^2}\right) \frac{ax}{\beta^2}\left(\frac{y}{ax}\right)^{\delta/2}$$

  and $I_{\delta/2-1}$ is the Bessel function with index $\delta/2 - 1$. This process (with here $a = 0.5$, $\beta = 3$, $\delta = 3$) is denoted by $\sqrt{\text{CIR}}$ since its square is actually a Cox-Ingersoll-Ross process. The estimation domain for this process is $[2, 10]^2$.

- A Cox-Ingersoll-Ross process, which is exactly the square of the previous process. The invariant distribution is a Gamma density with scale parameter $l = (1 - a^2)/(2\beta^2)$ and shape parameter $a = \delta/2$. The transition density is

$$\pi(x, y) = \frac{1}{2\beta^2} \exp\left(-\frac{y + a^2x}{2\beta^2}\right) I_{\delta/2-1}\left(\frac{a\sqrt{xy}}{\beta^2}\right)\left(\frac{y}{a^2x}\right)^{\delta/4-1/2}$$

  The used parameters are the following:

  (iii) $a = 3/4$, $\beta = \sqrt{7/48}$ and $\delta = 4$, estimated on $[0.1, 3]^2$.

(iv) $a = 1/3$, $\beta = 3/4$ and $\delta = 2$. This chain is estimated on $[0, 2]^2$.

- An ARCH process defined by $X_{n+1} = \sin(X_n) + (\cos(X_n) + 3)\varepsilon_{n+1}$ where the $\varepsilon_{n+1}$ are i.i.d. standard Gaussian. The transition density of this chain is
$$\pi(x, y) = \varphi\left(\frac{y - \sin(x)}{\cos(x) + 3}\right)\frac{1}{\cos(x) + 3}$$
and we estimate this process on $[-5, 5]^2$.

For this last chain, the stationary density is not explicit. So we simulate $n+500$ variables and we estimate only from the last $n$ to ensure the stationarity of the process. For the other chains, it is sufficient to simulate an initial variable $X_0$ with density $f$.

| $n$ | 50 | 100 | 250 | 500 | 1000 | noise |
|---|---|---|---|---|---|---|
| AR(i) | 0.579 | 0.407 | 0.270 | 0.230 | 0.209 | Lapl |
| | 0.599 | 0.480 | 0.313 | 0.272 | 0.245 | Gauss |
| AR(ii) | 0.389 | 0.294 | 0.195 | 0.155 | 0.139 | Lapl |
| | 0.339 | 0.304 | 0.280 | 0.273 | 0.271 | Gauss |
| $\sqrt{\text{CIR}}$ | 0.171 | 0.138 | 0.123 | 0.118 | 0.111 | Lapl |
| | 0.199 | 0.169 | 0.150 | 0.142 | 0.139 | Gauss |
| CIR(iii) | 0.420 | 0.345 | 0.237 | 0.195 | 0.175 | Lapl |
| | 0.337 | 0.302 | 0.276 | 0.245 | 0.209 | Gauss |
| CIR(iv) | 0.525 | 0.403 | 0.337 | 0.304 | 0.292 | Lapl |
| | 0.369 | 0.345 | 0.344 | 0.327 | 0.321 | Gauss |
| ARCH | 0.312 | 0.287 | 0.261 | 0.185 | 0.150 | Lapl |
| | 0.337 | 0.319 | 0.296 | 0.290 | 0.183 | Gauss |

TABLE 1

*MISE* $\mathbb{E}\|\Pi - \tilde{\Pi}\|^2$ *averaged over* $N = 200$ *samples.*

We consider two different noises:

**Laplace noise** In this case, the density of $\varepsilon_i$ is given by
$$q(x) = \frac{\lambda}{2}e^{-\lambda|x|}; \quad q^*(x) = \frac{\lambda^2}{\lambda^2 + x^2}; \quad \lambda = 5.$$

The smoothness parameter is $\gamma = 2$ so that the theoretical penalty is
$$\text{pen}(m) = C\Phi_1^2\|q\|_\infty^2 f_0^{-1}\frac{D_m^{10}}{n} = \frac{1}{n}\left(\frac{\lambda}{2}\right)^2 C\Phi_1^2 f_0^{-1}D_m^{10}$$

Several simulations lead to fix a constant $C$ very low. As the term $f_0^{-1}$ does not vary very much with regard to $C$, we choose to use the same following penalty for all the examples:
$$\text{pen}(m) = \frac{1}{n}\left(\frac{\lambda}{2}\right)^2\left(\frac{D_m}{4}\right)^{10}.$$

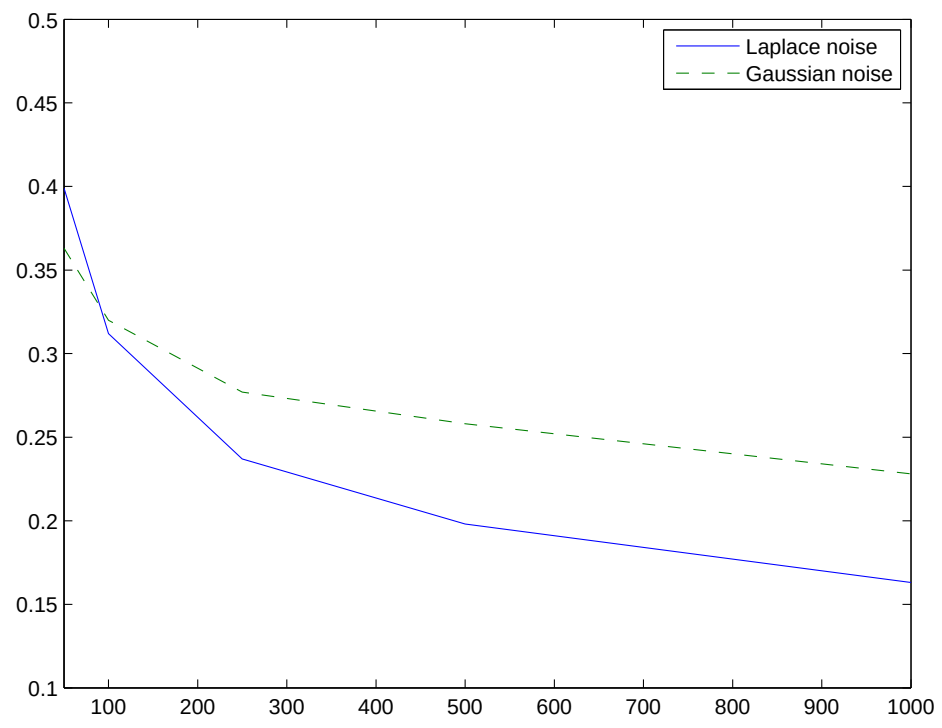

Fig 1. *Mean of the MISE for the six processes when n increases*

**Gaussian noise** In this case, the density of $\varepsilon_i$ is given by

$$q(x) = \frac{1}{\lambda\sqrt{2\pi}} e^{-\frac{x^2}{2\lambda^2}}; \quad q^*(x) = e^{-\frac{\lambda^2 x^2}{2}}; \quad \lambda = 0.3.$$

This noise does not verify Assumption H1 but it is interesting to see if this assumption is also necessary for practical purposes. Given the exponential regularity of this noise, we consider the following penalty

$$\mathrm{pen}(m) = \frac{\kappa}{n} \exp(\lambda^2 D_m^2)$$

where, by simulation experiments, we calibrate the penalty with $\kappa = 5$.

Table 1 presents the $L^2$ risk of our estimator of the transition density for the six Markov chains and the two noises. These results can be compared with those of Lacour (2007a) (Table 2) who studies the processes AR(i), $\sqrt{\mathrm{CIR}}$ and ARCH but directly observed, i.e. without noise. The risk values are then higher in our case, but with the same order, which is satisfactory. It is noticeable that the estimation works almost in the same way with the Gaussian noise, but with a slower decrease of the risk, as can be observed in Figure 1 . It is a classical phenomenon in deconvolution problems, since the Gaussian noise is much more regular than the Laplace noise.

Figure 2 allows visualizing the result for process ARCH observed through a Laplace noise: the surfaces $z = \Pi(x, y)$ and $z = \tilde{\Pi}(x, y)$ are presented. We also give figures of cross-sections of this kind of surfaces. We can see in Figure 3 the curves $z = \Pi(x, -0.44)$ versus $z = \tilde{\Pi}(x, -0.44)$ and the curves $z = \Pi(1.12, y)$ versus $z = \tilde{\Pi}(1.12, y)$ for the process AR(i). Generally, for a multidimensional estimation, the mixed control of the directions does not enable to do as well as a classical one-dimensional function estimation. Nevertheless here the curves are very close.

From a practical point of vue, it is difficult to compare the method described here and the one of Lacour (2007b). Indeed, the bases used are very different. However, we can say that the quotient method seems to give better results when the noise distribution is Gaussian (that is conform to theory). Nevertheless, the least squares procedure is better for a Laplace noise, especially when $n$ is small.

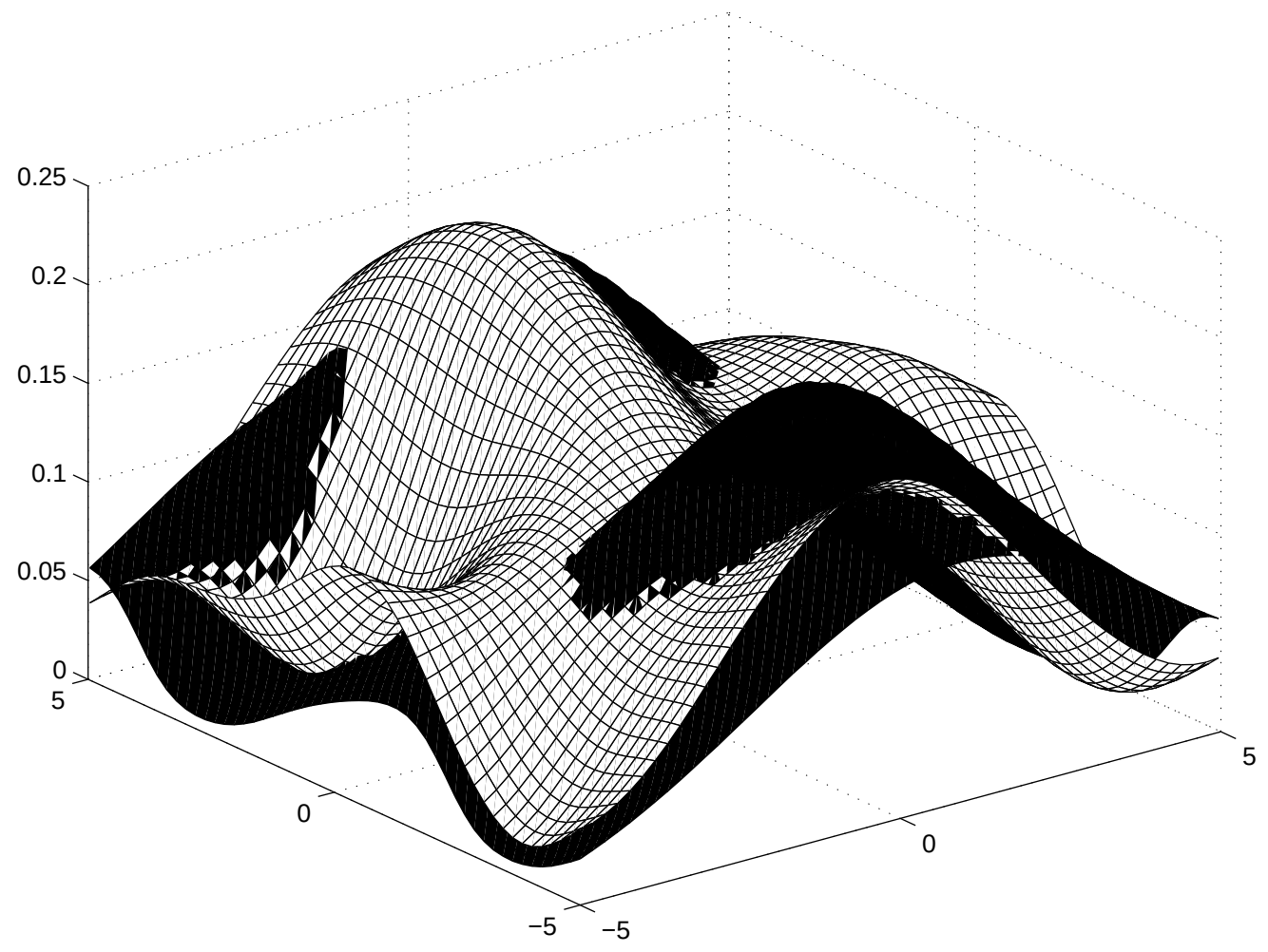

FIG 2. *True $\Pi$ (black) and estimator $\tilde{\Pi}$ (white) for process ARCH observed through a Laplace noise, $n = 500$*

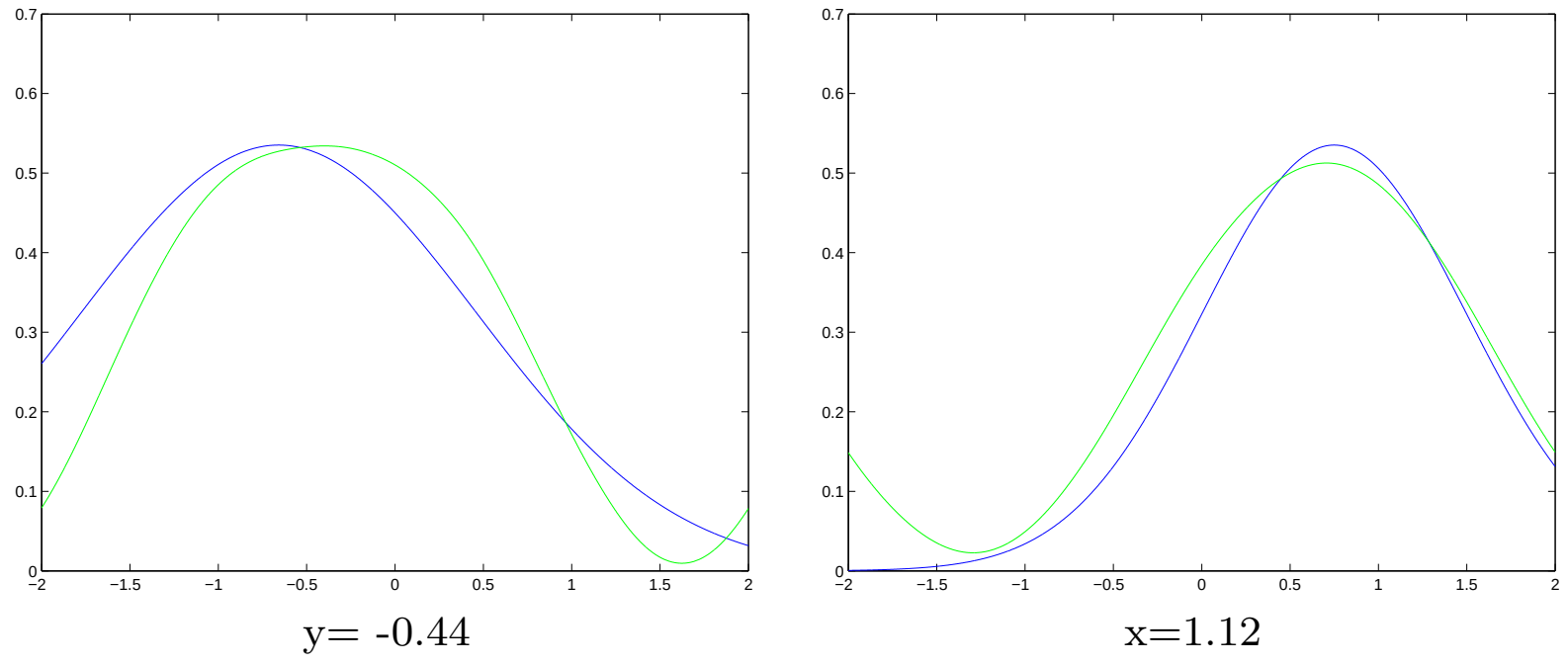

FIG 3. *Sections for process AR(i) observed through a Laplace noise, $n = 500$*

## 6. Detailed proofs

### *6.1. Proof of Lemma 1*

• Using

$$|\sum_k \varphi_{jk}^2(x)| \leq C(\varphi)2^{j/2}\|\varphi_{jk}\|_\infty \leq C'(\varphi)2^j, \tag{8}$$

P3 holds if $\Phi_1 \geq 2C'(\varphi)$.

• The computation of $v_{\varphi_{jk}}$ gives

$$|v_{\varphi_{jk}}(x)| \leq \frac{2^{j/2}}{2\pi} \int \frac{|\varphi^*(v)|}{|q^*(-v2^j)|} dv$$

Next, it follows from Assumption H1 that $|v_{\varphi_{jk}}(x)| \leq C_{1,\gamma}(2^j)^{\gamma+1/2}/2\pi k_0$ using Lemma 5 (Section 6.8) since $r > \gamma + 1$. Then, for all $x$, $\sum_k |v_{\varphi_{jk}}(x)|^2 \leq 3.2^j C_{1,\gamma}^2 \frac{k_0^{-2}}{4\pi^2}(2^j)^{2\gamma+1}$ that establishes P4 with $\Phi_1 \geq 3C_{1,\gamma}^2 k_0^{-2}/(4\pi^2)$.

• To prove P5, we apply the Parseval equality. That yields

$$\int |v_{\varphi_{jk}}|^2 = \frac{1}{2\pi} \int \frac{|\varphi^*(v)|^2}{|q^*(-v2^j)|^2} dv.$$

Using H1 and given that $2r > 2\gamma + 1$, we obtain $\int |v_{\varphi_{jk}}|^2 \leq C_{2,2\gamma}(2^j)^{2\gamma}/2\pi k_0^2$ And finally P5 holds with $\Phi_1 \geq 3C_{2,2\gamma} k_0^{-2}/(2\pi)$.

• We begin with computing $|v_{\varphi_{jk}\varphi_{jk'}}(x)|$ by using that $(\varphi_{jk}\varphi_{jk'})^*$ is equal to the convolution product $\varphi_{jk}^* * \varphi_{jk'}^*$.

$$\begin{aligned} |v_{\varphi_{jk}\varphi_{jk'}}(x)| &\leq \frac{2^{-j}}{2\pi} \iint \frac{|\varphi^*(v/2^j)||\varphi^*((u-v)/2^j)|}{|q^*(-u)|} dudv \\ &\leq \frac{k_0^{-1}}{2\pi} 2^j (2^j)^\gamma \iint |\varphi^*(y)\varphi^*(x-y)|(x^2+1)^{\gamma/2} dxdy. \end{aligned}$$

Then Lemma 5 (Section 6.8) shows that

$$\begin{aligned} |v_{\varphi_{jk}\varphi_{jk'}}(x)| \quad &\leq \quad \frac{k_0^{-1}}{2\pi}(2^j)^{\gamma+1} C_r \Big[ \int_{|x|>1} |x|^{1-r}(x^2+1)^{\gamma/2} dx \\ &\quad + \int_{|x|\leq 1} (x^2+1)^{\gamma/2} dx \Big]. \end{aligned}$$

Hence, since $r > \gamma + 2$, there exists $C > 0$ such that $|v_{\varphi_{jk}\varphi_{jk'}}(x)| \leq C(2^j)^{\gamma+1}$. The fact that $\varphi_{jk}$ and $\varphi_{jk'}$ have disjoint supports if $k + N \leq k' - N + 1$ or $k' + N \leq k - N + 1$ enables to prove P6 with $\Phi_1 \geq 3(4N-3)C^2$.

• Applying Parseval's equality,

$$\int |v_{\varphi_{jk}\varphi_{jk'}}|^2 = \frac{2^j}{2\pi} \int \frac{|(\varphi_{jk}\varphi_{jk'})^*|^2(2^j v)}{|q^*(-2^j v)|^2} dv.$$

But,using Lemma 5,

$$|(\varphi_{jk}\varphi_{jk'})^*(2^j v)| \leq \int |\varphi^*(y)||\varphi^*(v-y)|dy \leq C_r \left[|v|^{1-r}\mathbb{1}_{|v|>1} + \mathbb{1}_{|v|\leq 1}\right] \quad (9)$$

Then, it follows that

$$\int |v_{\varphi_{jk}\varphi_{jk'}}|^2 \leq \frac{k_0^{-2}2^j}{2\pi} C_r^2 \int (|v|^{2(1-r)}\mathbb{1}_{|v|>1} + \mathbb{1}_{|v|\leq 1})((2^j v)^2+1)^\gamma dv \leq C(2^j)^{2\gamma+1}.$$

It is then sufficient to sum this quantity for all $k, k'$ by taking into account the superposition of the supports to prove P7 as soon as $\Phi_1 \geq 3C(4N-3)$.

### *6.2. Proof of Lemma 2*

1. First we write that

$$V_T(Y_k, Y_{k+1}) = \frac{1}{4\pi^2}\int e^{iY_k u + iY_{k+1}v}\frac{T^*(u,v)}{q^*(-u)q^*(-v)}dudv$$

so that, by denoting $\mathbb{X} = (X_1, \dots, X_{n+1})$,

$$\mathbb{E}[V_T(Y_k, Y_{k+1})|\mathbb{X}] = \frac{1}{4\pi^2}\int \mathbb{E}[e^{iY_k u + iY_{k+1}v}|\mathbb{X}]\frac{T^*(u,v)}{q^*(-u)q^*(-v)}dudv.$$

By using the independence between $(X_i)$ and $(\varepsilon_i)$, we compute

$$\begin{aligned}\mathbb{E}[e^{iY_k u + iY_{k+1}v}|\mathbb{X}] = \mathbb{E}[e^{iX_k u + iX_{k+1}v}e^{i\varepsilon_k u}e^{i\varepsilon_{k+1}u}|\mathbb{X}]\\ = e^{iX_k u + iX_{k+1}v}\mathbb{E}[e^{i\varepsilon_k u}]\mathbb{E}[e^{i\varepsilon_{k+1}v}] = e^{iX_k u + iX_{k+1}v}q^*(-u)q^*(-v).\end{aligned}$$

Then

$$\mathbb{E}[V_T(Y_k, Y_{k+1})|\mathbb{X}] = \frac{1}{4\pi^2}\int e^{iX_k u + iX_{k+1}v}T^*(u,v)dudv = T(X_k, X_{k+1}).$$

3. We proceed in a similar way for $Q$. Since $Q_T(Y_k) = (1/2\pi)\int e^{iY_k u}T^*(u,0)$ $(q^*(-u))^{-1}du$, then

$$\mathbb{E}[Q_T(Y_k)|\mathbb{X}] = \frac{1}{2\pi}\int \mathbb{E}[e^{iY_k u}|\mathbb{X}]\frac{T^*(u,0)}{q^*(-u)}du.$$

By using the independence between $(X_i)$ and $(\varepsilon_i)$, we compute

$$\mathbb{E}[e^{iY_k u}|\mathbb{X}] = \mathbb{E}[e^{iX_k u}e^{i\varepsilon_k u}|\mathbb{X}] = e^{iX_k u}\mathbb{E}[e^{i\varepsilon_k u}] = e^{iX_k u}q^*(-u).$$

Thus

$$\mathbb{E}[Q_T(Y_k)|\mathbb{X}] = \frac{1}{2\pi}\int e^{iX_k u}q^*(-u)\frac{T^*(u,0)}{q^*(-u)}du = \frac{1}{2\pi}\int e^{iX_k u}T^*(u,0)du.$$

By denoting by $T_y$ the function $x \mapsto T_y(x) = T(x,y)$, we obtain

$$T^*(u,0) = \iint e^{-ixu}T_y(x)dxdy = \int T_y^*(u)dy$$

and then

$$\begin{aligned}\frac{1}{2\pi}\int e^{iX_k u}T^*(u,0)du &= \frac{1}{2\pi}\iint e^{iX_k u}T_y^*(u)dydu\\ &= \int T_y(X_k)dy = \int T(X_k, y)dy.\end{aligned} \tag{10}$$

### 6.3. Proof of Theorem 1

We start with introducing some auxiliary variables whose existence is ensured by Assumption H5 of mixing. In the case of arithmetical mixing, since $\theta > 8$, there exists a real $c$ such that $0 < c < 3/8$ and $c\theta > 3$. In this case, we set $q_n = \lfloor n^c \rfloor$. In the case of geometric mixing, we set $q_n = \lfloor c\ln(n) \rfloor$ where $c$ is a real larger than $3/\theta$.

For the sake of simplicity, we suppose that $n+1 = 2p_nq_n$, with $p_n$ an integer. Let for $l = 0,\dots,p_n-1$, $A_l = (X_{2lq_n+1},...,X_{(2l+1)q_n})$ and $B_l = (X_{(2l+1)q_n+1},...,X_{(2l+2)q_n})$. As in Viennet (1997), by using Berbee's coupling Lemma, we can build a sequence $(A_l^*)$ such that

$$\begin{cases} A_l \text{ and } A_l^* \text{ have the same distribution,} \\ A_l^* \text{ and } A_{l'}^* \text{ are independent if } l \neq l', \\ P(A_l \neq A_l^*) \leq \beta_{q_n}. \end{cases} \tag{11}$$

In the same way, we build $(B_l^*)$ and we define for any $l \in \{0,\dots,p_n-1\}$, $A_l^* = (X^*_{2lq_n+1},...,X^*_{(2l+1)q_n})$, $B_l^* = (X^*_{(2l+1)q_n+1},...,X^*_{(2l+2)q_n})$ so that the sequence $(X_1^*,\dots,X_n^*)$ is well defined. We can now define

$$\Omega_X^* = \{\forall i, 1 \leq i \leq n+1 \quad X_i = X_i^*\}.$$

Let us recall that $\mathcal{S}$ is the space $\mathbb{S}_m$ with maximal dimension $\mathcal{D}^2 \leq n^{\frac{1}{4\gamma+2}}$. We now adopt the notations

$$\Delta = \{\forall T \in \mathcal{S} \quad \|T\|_f^2 \leq \frac{3}{2}\Psi_n(T)\}; \qquad \Omega = \Delta \cap \Omega_X^*.$$

Let us fix $m \in \mathcal{M}_n$. We denote by $\Pi_m$ the orthogonal projection of $\Pi$ on $\mathbb{S}_m$. Then we have the decomposition

$$\begin{aligned}\mathbb{E}\|\tilde{\Pi}-\Pi\|_A^2 \quad \leq \quad & 2\mathbb{E}\left(\|\tilde{\Pi}-\Pi_m\|_A^2 \mathbb{1}_\Omega \mathbb{1}_{\|\hat{\Pi}_{\hat{m}}\| \leq k_n}\right) + 2\mathbb{E}\left(\|\tilde{\Pi}-\Pi_m\|_A^2 \mathbb{1}_\Omega \mathbb{1}_{\|\hat{\Pi}_{\hat{m}}\| > k_n}\right) \\ & + 2\mathbb{E}\left(\|\tilde{\Pi}-\Pi_m\|_A^2 \mathbb{1}_{\Omega^c}\right) + 2\|\Pi_m - \Pi\|_A^2 \\ \leq \quad & 2\mathbb{E}\left(\|\hat{\Pi}_{\hat{m}}-\Pi_m\|_A^2 \mathbb{1}_\Omega\right) + 2\|\Pi_m\|_A^2 \mathbb{E}\left(\mathbb{1}_\Omega \mathbb{1}_{\|\hat{\Pi}_{\hat{m}}\| > k_n}\right) \\ & + 2\mathbb{E}\left([2\|\tilde{\Pi}\|_A^2 + 2\|\Pi_m\|_A^2]\mathbb{1}_{\Omega^c}\right) + 2\|\Pi_m - \Pi\|_A^2.\end{aligned}$$

Now, using the Markov inequality and the definition of $\tilde{\Pi}$,

$$\begin{aligned}\mathbb{E}\|\tilde{\Pi}-\Pi\|_A^2 \leq 2\mathbb{E}\left(\|\hat{\Pi}_{\hat{m}}-\Pi_m\|_A^2 \mathbb{1}_\Omega\right) + 2\|\Pi\|_A^2 \frac{\mathbb{E}(\|\hat{\Pi}_{\hat{m}}\|^2 \mathbb{1}_\Omega)}{k_n^2} \\ + 4(k_n^2 + \|\Pi\|_A^2)\mathbb{E}\left(\mathbb{1}_{\Omega^c}\right) + 2\|\Pi_m - \Pi\|_A^2.\end{aligned}$$

But $\mathbb{E}(\|\hat{\Pi}_{\hat{m}}\|^2 \mathbb{1}_\Omega) \leq 2\mathbb{E}(\|\hat{\Pi}_{\hat{m}}-\Pi_m\|_A^2 \mathbb{1}_\Omega) + 2\|\Pi_m\|_A^2$ and $k_n = \sqrt{n}$, so

$$\begin{aligned}\mathbb{E}\|\tilde{\Pi}-\Pi\|_A^2 \leq 2\mathbb{E}\left(\|\hat{\Pi}_{\hat{m}}-\Pi_m\|_A^2 \mathbb{1}_\Omega\right)(1 + 2\|\Pi\|_A^2) + \frac{4\|\Pi\|_A^4}{n} \\ + 4(n + \|\Pi\|_A^2)P(\Omega^c) + 2\|\Pi_m - \Pi\|_A^2.\end{aligned}$$

We now state the following proposition :

**Proposition 2.** *There exists $C_0 > 0$ such that*

$$P(\Omega^c) \leq \frac{C_0}{n^2}.$$

Hence

$$\begin{aligned}\mathbb{E}\|\tilde{\Pi} - \Pi\|_A^2 \leq 2\|\Pi_m - \Pi\|_A^2 + 2\mathbb{E}\left(\|\hat{\Pi}_{\hat{m}} - \Pi_m\|_A^2 \mathbb{1}_\Omega\right)(1 + 2\|\Pi\|_A^2) \\ + \frac{4}{n}(\|\Pi\|_A^4 + C_0(1 + \|\Pi\|_A^2)).\end{aligned} \tag{12}$$

Now we have to bound $\mathbb{E}\left(\|\hat{\Pi}_{\hat{m}} - \Pi_m\|_A^2 \mathbb{1}_\Omega\right)$. The estimators $\hat{\Pi}_m$ are defined by minimization of the contrast on a set $\Gamma$ defined in (6). Let us prove that this set $\Gamma$ contains $\Omega$. More precisely, we prove that $\Delta \subset \Gamma$. For $T = \sum_\lambda a_\lambda \omega_\lambda \in \mathbb{S}_m$, the matrix $A_m = (a_\lambda)$ of its coefficients in the basis $(\omega_\lambda(x, y))$ verifies $\Psi_n(T) = {}^t A_m G_m A_m$. Then, on $\Delta$,

$${}^t A_m G_m A_m \geq \frac{2}{3}\|T\|_f^2 \geq \frac{2}{3} f_0 \|T\|^2.$$

Now, using P2, $\|T\|^2 = {}^t A_m A_m$ and then ${}^t A_m G_m A_m \geq (2/3) f_0 {}^t A_m A_m$. If $\mu$ is an eigenvalue of $G_m$, there exists $A_m \neq 0$ such that $G_m A_m = \mu A_m$ and then ${}^t A_m G_m A_m = \mu {}^t A_m A_m$. Then, on $\Delta$,

$$\mu {}^t A_m A_m \geq \frac{2}{3} f_0 {}^t A_m A_m.$$

Consequently $\mu \geq (2/3) f_0$. So $\Delta \subset \Gamma$ and $\hat{\Pi}_{\hat{m}}$ minimizes the contrast on $\Delta$.

We now observe that, for all functions $T, S$

$$\gamma_n(T) - \gamma_n(S) = \Psi_n(T - S) - \frac{2}{n}\sum_{k=1}^{n}[V_{(T-S)}(Y_k, Y_{k+1}) - Q_{(T-S)S}(Y_k)].$$

Then, since on $\Delta$, $\gamma_n(\hat{\Pi}_{\hat{m}}) + \mathrm{pen}(\hat{m}) \leq \gamma_n(\Pi_m) + \mathrm{pen}(m)$,

$$\begin{aligned}\Psi_n(\hat{\Pi}_{\hat{m}} - \Pi_m) &\leq \frac{2}{n}\sum_{k=1}^{n}[V_{(\hat{\Pi}_{\hat{m}} - \Pi_m)}(Y_k, Y_{k+1}) - Q_{(\hat{\Pi}_{\hat{m}} - \Pi_m)\Pi_m}(Y_k)] \\ &\quad + \mathrm{pen}(m) - \mathrm{pen}(\hat{m}) \\ &\leq 2Z_{n,m}(\hat{\Pi}_{\hat{m}} - \Pi_m) + \mathrm{pen}(m) - \mathrm{pen}(\hat{m}) \\ &\leq 2\|\hat{\Pi}_{\hat{m}} - \Pi_m\|_f \sup_{T \in B_f(m, \hat{m})} Z_{n,m}(T) + \mathrm{pen}(m) - \mathrm{pen}(\hat{m})\end{aligned}$$

where

$$Z_{n,m}(T) = \frac{1}{n}\sum_{k=1}^{n}[V_T(Y_k, Y_{k+1}) - Q_{T\Pi_m}(Y_k)]$$

and, for all $m'$, $B_f(m,m') = \{T \in S_m + S_{m'}, \quad \|T\|_f = 1\}$. Now let $p(.,.)$ be a function such that for all $m, m'$, $12p(m,m') \leq \text{pen}(m) + \text{pen}(m')$. Then

$$\Psi_n(\hat{\Pi}_{\hat{m}} - \Pi_m) \leq \frac{1}{3}\|\hat{\Pi}_{\hat{m}} - \Pi_m\|_f^2 + 3[\sup_{T \in B_f(m,\hat{m})} Z_{n,m}^2(T) - 4p(m,\hat{m})] + 2\text{pen}(m).$$

So, using the definition of $\Delta \supset \Omega$,

$$\|\hat{\Pi}_{\hat{m}} - \Pi_m\|_f^2 \mathbb{1}_\Omega \leq \frac{3}{2}\Psi_n(\hat{\Pi}_{\hat{m}} - \Pi_m)\mathbb{1}_\Omega$$
$$\leq \frac{1}{2}\|\hat{\Pi}_{\hat{m}} - \Pi_m\|_f^2 \mathbb{1}_\Omega + \frac{9}{2}\sum_{m' \in \mathcal{M}_n}[\sup_{T \in B_f(m,m')} Z_{n,m}^2(T) - 4p(m,m')]\mathbb{1}_\Omega + 3\text{pen}(m)$$

Thus

$$\frac{1}{2}\|\hat{\Pi}_{\hat{m}} - \Pi_m\|_f^2 \mathbb{1}_\Omega \leq \frac{9}{2}\sum_{m' \in \mathcal{M}_n}[\sup_{T \in B_f(m,m')} Z_{n,m}^2(T) - 4p(m,m')]\mathbb{1}_\Omega + 3\text{pen}(m)$$

And using Assumption H3,

$$\|\hat{\Pi}_{\hat{m}} - \Pi_m\|_A^2 \mathbb{1}_\Omega \leq 9f_0^{-1}\sum_{m' \in \mathcal{M}_n}[\sup_{T \in B_f(m,m')} Z_{n,m}^2(T) - 4p(m,m')]\mathbb{1}_\Omega + 6f_0^{-1}\text{pen}(m) \tag{13}$$

Now, by denoting $\mathbb{E}_X$ the expectation conditionally to $X_1, \ldots, X_{n+1}$, the process $Z_{n,m}(T)$ can be split in the following way :

$$Z_{n,m}(T) = Z_n^{(1)}(T) - Z_n^{(2)}(T) + Z_n^{(3)}(T) + \iint T(x,y)(\Pi - \Pi_m)(x,y)f(x)dxdy$$

with

$$\begin{cases} Z_n^{(1)}(T) = \frac{1}{n}\sum_{k=1}^n V_T(Y_k, Y_{k+1}) - \mathbb{E}_X[V_T(Y_k, Y_{k+1})] \\ Z_n^{(2)}(T) = \frac{1}{n}\sum_{k=1}^n Q_{T\Pi_m}(Y_k) - \mathbb{E}[Q_{T\Pi_m}(Y_k)] \\ Z_n^{(3)}(T) = \frac{1}{n}\sum_{k=1}^n T(X_k, X_{k+1}) - \mathbb{E}[T(X_k, X_{k+1})] \end{cases}$$

Then, by introducing functions $p_1(.,.)$, $p_2(.,.)$ and $p_3(.,.)$

$$\sup_{T \in B_f(m,m')} Z_{n,m}^2(T) - 4p(m,m') \leq 4\sup_{T \in B_f(m,m')}(Z_n^{(1)}(T)^2 - p_1(m,m'))$$
$$+4\sup_{T \in B_f(m,m')}(Z_n^{(2)}(T)^2 - p_2(m,m')) + 4\sup_{T \in B_f(m,m')}(Z_n^{(3)}(T)^2 - p_3(m,m'))$$
$$+4((p_1 + p_2 + p_3)(m,m') - p(m,m')) + 4\sup_{T \in B_f(m,m')}\|(\Pi - \Pi_m)\mathbb{1}_A\|_f^2\|T\|_f^2$$

We now use the following propositions.

**Proposition 3.** *Let* $p_1(m,m') = K_1(\gamma)\Phi_1^2 f_0^{-1}\|q\|_\infty^2 D_{m''}^{4\gamma+2}/n$ *where* $m'' = \max(m,m')$. *Then, if* $r > 2\gamma + 1/2$, *there exists a positive constant* $C_1$ *such that*

$$\sum_{m'\in\mathcal{M}_n} \mathbb{E}\left(\left[\sup_{T\in B_f(m,m')} Z_n^{(1)}(T)^2 - p_1(m,m')\right]_+\right) \le \frac{C_1}{n}.$$

**Proposition 4.** *Let* $p_2(m,m') = p_2^{(1)}(m,m') + p_2^{(2)}(m,m')$ *with* $p_2^{(1)}(m,m') = K_2\|\Pi\|_A^2\ D_{m''}^{2\gamma+7/2}/n$ *and* $p_2^{(2)}(m,m') = K_2\|\Pi\|_A^2(\sum_k \beta_k)D_{m''}^3/n$ *where* $m'' = \max(m,m')$. *Then, if* $r > 2\gamma+3/2$, *there exists a positive constant* $C_2$ *such that*

$$\sum_{m'\in\mathcal{M}_n} \mathbb{E}\left(\left[\sup_{T\in B_f(m,m')} Z_n^{(2)}(T)^2 - p_2(m,m')\right]_+ \mathbb{1}_\Omega\right) \le \frac{C_2}{n}.$$

**Proposition 5.** *Let* $p_3(m,m') = K_3\sum_k \beta_{2k} D_{m''}^2/n$ *where* $m'' = \max(m,m')$. *Then, there exists a positive constant* $C_3$ *such that*

$$\sum_{m'\in\mathcal{M}_n} \mathbb{E}\left(\left[\sup_{T\in B_f(m,m')} Z_n^{(3)}(T)^2 - p_3(m,m')\right]_+ \mathbb{1}_\Omega\right) \le \frac{C_3}{n}.$$

The first two propositions are proved in Sections 6.6 and 6.7. The last proposition is proved in Lacour (2007b) Section 6.5 (for another basis but only the property P3 $\|\sum_{jk}\varphi_{jk}^2\|_\infty \le \Phi_1 D_m$ is used).

Then we get

$$\sum_{m'\in\mathcal{M}_n} \mathbb{E}\left(\left[\sup_{T\in B_f(m,m')} Z_{n,m}^2(T) - 4p(m,m')\right]\mathbb{1}_\Omega\right) \le 4\frac{C_1+C_2+C_3}{n}$$
$$+4\|(\Pi-\Pi_m)\mathbb{1}_A\|_f^2 + 4\sum_{m'\in\mathcal{M}_n}((p_1+p_2+p_3)(m,m') - p(m,m')).$$

But, if $\gamma > 3/4$, $4\gamma+2 > 2\gamma+7/2$ and there exists $m_2$ such that for all $m' > m_2$, $p_1(m,m') > p_2(m,m') + p_3(m,m')$. It implies that

$$\sum_{m'\in\mathcal{M}_n}(p_1(m,m') + p_2(m,m') + p_3(m,m') - 2p_1(m,m'))$$
$$\le \sum_{m'\le m_2}(p_2(m,m') + p_3(m,m') - p_1(m,m')) \le \frac{C(m_2)}{n}.$$

Thus in the case $\gamma > 3/4$, we choose $p = 2p_1$ and

$$\sum_{m'\in\mathcal{M}_n} \mathbb{E}\left(\left[\sup_{T\in B_f(m')} Z_{n,m}^2(T) - 4p(m,m')\right]\mathbb{1}_\Omega\right) \le 4\frac{C_1+C_2+C_3+C(m_2)}{n}$$
$$+4\|f\|_{\infty,A_1}\|\Pi-\Pi_m\|_A^2 \tag{14}$$

If $\gamma = 3/4$, we choose $p = 2(p_1 + p_2^{(1)})$. Since there exists $m_2$ such that for all $m' > m_2$, $p_1(m, m') + p_2^{(1)}(m, m') > p_2^{(2)}(m, m') + p_3(m, m')$, we can write

$$\sum_{m' \in \mathcal{M}_n} (p_1(m, m') + p_2(m, m') + p_3(m, m') - p(m, m'))$$
$$\leq \sum_{m' \leq m_2} (p_2^{(2)}(m, m') + p_3(m, m') - p_1(m, m') - p_2^{(1)}(m, m')) \leq \frac{C(m_2)}{n}$$

and (14) holds.

Finally, combining (12), (13) and (14), we obtain

$$\mathbb{E}\|\tilde{\Pi} - \Pi\|_A^2 \leq 2\|\Pi_m - \Pi\|_A^2 + \frac{4}{n}(\|\Pi\|_A^4 + C_0(1 + \|\Pi\|_A^2))$$
$$+2(1 + 2\|\Pi\|_A^2)9f_0^{-1}\Big[4\frac{C_1 + C_2 + C_3 + C(m_2)}{n} + 4\|f\|_{\infty,A_1}\|\Pi - \Pi_m\|_A^2\Big]$$
$$+2(1 + 2\|\Pi\|_A^2)6f_0^{-1}\text{pen}(m).$$

Then, by letting $C = \max(2 + 72f_0^{-1}\|f\|_{\infty,A_1}(1 + 2\|\Pi\|_A^2), 12f_0^{-1}(1 + 2\|\Pi\|_A^2))$,

$$\mathbb{E}\|\tilde{\Pi} - \Pi\|_A^2 \leq C \inf_{m \in \mathcal{M}_n} (\|\Pi_m - \Pi\|_A^2 + \text{pen}(m)) + \frac{C'}{n}$$

We still have to verify that $12p(m, m') \leq \text{pen}(m) + \text{pen}(m')$. But, if $\gamma > 3/4$,

$$12p(m, m') = 24K_1\frac{D_{m''}^{4\gamma+2}}{n} = 24K_1\frac{\dim(S_m + S_{m'})^{4\gamma+2}}{n} \leq \text{pen}(m) + \text{pen}(m')$$

with $\text{pen}(m) \geq 24K_1 D_m^{4\gamma+2}/n$. And if $\gamma = 3/4$,

$$12p(m, m') = 24(K_1 + K_2\|\Pi\|_A^2)\frac{D_{m''}^5}{n} \leq \text{pen}(m) + \text{pen}(m')$$

with $\text{pen}(m) \geq 24(K_1 + K_2\|\Pi\|_A^2)D_m^{4\gamma+2}/n$.

### *6.4. Proof of Corollary 1*

It follows from Meyer (1990) Chapter 6, Section 10 that $\Pi$ belongs to $B_{2,\infty}^{\alpha}$ if and only if $\sup_{j \geq J} 2^{2j\alpha}(\sum_{k,l}|a_{jkl}|^2)^{1/2} < \infty$ with $a_{jkl} = \int \Pi(x, y)\varphi_{jk}(x)\varphi_{jl}(y)dxdy$. Then

$$d_A^2(\Pi, S_m) = \sum_{j > m}\sum_{k,l}|a_{jkl}|^2 \leq C\sum_{j > m} 2^{-4j\alpha} \leq C' D_m^{-2\alpha}$$

Since $d_A^2(\Pi, S_m) = O(D_m^{-2\alpha})$, Theorem 1 becomes

$$\mathbb{E}\|\tilde{\Pi} - \Pi\|_A^2 \leq C'' \inf_{m \in \mathcal{M}_n}\{D_m^{-2\alpha} + \frac{D_m^{4\gamma+2}}{n}\}.$$

with $C''$ a positive constant. By setting $D_{m_1}$ the integer part of $n^{1/(4\gamma+2\alpha+2)}$, then

$$\mathbb{E}\|\tilde{\Pi} - \Pi\|_A^2 \leq C''\{D_{m_1}^{-2\alpha} + \frac{D_{m_1}^{4\gamma+2}}{n}\} = O(n^{-\frac{2\alpha}{4\gamma+2\alpha+2}}).$$

### *6.5. Proof of Proposition 2*

We first remark that $P(\Omega^c) \leq P(\Omega_X^{*c}) + P(\Delta^c \cap \Omega_X^*)$. In the geometric case $\beta_{q_n} \leq e^{-\theta c \ln(n)} \leq n^{-\theta c}$ and in the other case $\beta_{q_n} \leq (q_n)^{-\theta} \leq n^{-\theta c}$. Then

$$P(\Omega_X^{*c}) \leq 2p_n\beta_{q_n} \leq n^{1-c\theta}.$$

But, $c\theta > 3$ and so $P(\Omega_X^{*c}) \leq n^{-2}$. We still have to bound $P(\Delta^c \cap \Omega_X^*)$. To do this, we observe that if $\omega \in \Delta^c$, then there exists $T$ in $\mathcal{S}$ such that $\|T\|_f^2 > (3/2)\Psi_n(T)$ and then $\|T\|_f^2 > (3/2)\mathbb{E}_X\Psi_n(T)$. But $\mathbb{E}_X\Psi_n(T) = \frac{1}{n}\sum_{k=1}^n \int T^2(X_k, y)dy$. So $P(\Delta^c \cap \Omega_X^*) \leq P(\Delta'^c \cap \Omega_X^*)$ with

$$\Delta' = \{\forall T \in \mathcal{S} \quad \|T\|_f^2 \leq \frac{3}{2}\frac{1}{n}\sum_{k=1}^n \int T^2(X_k, y)dy\}.$$

Let us remark that $(1/n)\sum_{k=1}^n \int T^2(X_k, y)dy - \|T\|_f^2 = \nu_n(T^2)$ with

$$\nu_n(T) = \frac{1}{n}\sum_{i=1}^n \int [T(X_i, y) - \mathbb{E}(T(X_i, y))]dy.$$

Hence

$$P(\Delta'^c \cap \Omega_X^*) \leq P(\sup_{T\in\mathcal{B}} |\nu_n(T^2)|\mathbb{1}_{\Omega_X^*} > 1/3)$$

with $\mathcal{B} = \{T \in \mathcal{S} \quad \|T\|_f = 1\}$.

A function $T$ in $\mathcal{S}$ can be written $T(x, y) = \sum_{j=J}^{m_0}\sum_{kl} a_{jkl}\varphi_{jk}(x)\varphi_{jl}(y)$ where $m_0$ is such that $\mathcal{S} = \mathbb{S}_{m_0}$. Then

$$\nu_n(T^2)\mathbb{1}_{\Omega_X^*} = \sum_{jkk'l} a_{jkl}a_{jk'l}\bar{\nu}_n(\varphi_{jk}\varphi_{jk'})$$

where

$$\bar{\nu}_n(u) = \frac{1}{n}\sum_{i=1}^n [u(X_i^*) - \mathbb{E}(u(X_i^*))]. \tag{15}$$

Let $b_{jk} = (\sum_l a_{jkl}^2)^{1/2}$, then $|\nu_n(T^2)|\mathbb{1}_{\Omega_X^*} \leq \sum_{jkk'} b_{jk}b_{jk'}|\bar{\nu}_n(\varphi_{jk}\varphi_{jk'})|$ and, if $T \in \mathcal{B}$, $\sum_{jk} b_{jk}^2 = \sum_{jkl} a_{jkl}^2 = \|T\|^2 \leq f_0^{-1}$

Thus,

$$\sup_{T\in\mathcal{B}} |\nu_n(T^2)|\mathbb{1}_{\Omega_X^*} \leq f_0^{-1} \sup_{\sum b_{jk}^2=1} \sum_{jkk'} b_{jk}b_{jk'}|\bar{\nu}_n(\varphi_{jk}\varphi_{jk'})|.$$

For the sake of simplicity, we denote $\lambda = (j, k)$ and $\lambda' = (j, k')$ so that

$$\sup_{T\in\mathcal{B}} |\nu_n(T^2)|\mathbb{1}_{\Omega_X^*} \leq f_0^{-1} \sup_{\sum b_\lambda^2=1} \sum_{\lambda\lambda'} b_\lambda b_{\lambda'}|\bar{\nu}_n(\varphi_\lambda\varphi_{\lambda'})|.$$

**Lemma 3.** *Let* $B_{\lambda,\lambda'} = \|\varphi_\lambda\varphi_{\lambda'}\|_\infty$ *and* $V_{\lambda,\lambda'} = \|\varphi_\lambda\varphi_{\lambda'}\|_2$. *Let, for any symmetric matrix* $(A_{\lambda,\lambda'})$

$$\bar\rho(A) = \sup_{\sum a_\lambda^2=1} \sum_{\lambda,\lambda'} |a_\lambda a_{\lambda'}| A_{\lambda,\lambda'}$$

*and* $L(\varphi) = \max\{\bar\rho^2(V), \bar\rho(B)\}$. *Then there exists* $\Phi_0 > 0$ *such that* $L(\varphi) \leq \Phi_0\mathcal{D}^2$.

This lemma is proved in Baraud et al. (2001) for an orthonormal basis verifying $\|\sum_\lambda \varphi_\lambda^2\|_\infty \leq \Phi_0\mathcal{D}$, that is ensured by property P3.

Now let $x = \dfrac{f_0^2}{24\|f\|_{\infty,A_1}L(\varphi)}$ and

$$D = \left\{\forall\lambda\forall\lambda' \quad |\bar\nu_n(\varphi_\lambda\varphi_{\lambda'})| \leq \left[B_{\lambda,\lambda'}x + V_{\lambda,\lambda'}\sqrt{2\|f\|_{\infty,A_1}x}\right]\right\}.$$

On $D$:

$$\sup_{T\in\mathcal{B}} |\nu_n(T^2)|\mathbb{1}_{\Omega_X^*} \leq f_0^{-1} \sup_{\sum b_\lambda^2=1} \sum_{\lambda,\lambda'} b_\lambda b_{\lambda'} \left[B_{\lambda,\lambda'}x + V_{\lambda,\lambda'}\sqrt{2\|f\|_{\infty,A_1}x}\right]$$

$$\leq f_0^{-1}\left[\bar\rho(B)x + \bar\rho(V)\sqrt{2\|f\|_{\infty,A_1}x}\right]$$

$$\leq \frac{f_0}{24\|f\|_{\infty,A_1}}\frac{\bar\rho(B)}{L(\varphi)} + \frac{1}{\sqrt{12}}\left(\frac{\bar\rho^2(V)}{L(\varphi)}\right)^{1/2} \leq \frac{1}{24} + \frac{1}{2\sqrt3} < \frac13.$$

Then $P\left(\sup_{T\in\mathcal{B}} |\nu_n(T^2)|\mathbb{1}_{\Omega_X^*} > 1/3\right) \leq P(D^c)$. But $\bar\nu_n(u) = \bar\nu_{n,1}(u)/2 + \bar\nu_{n,2}(u)/2$ with

$$\bar\nu_{n,s}(u) = \frac{1}{p_n}\sum_{l=0}^{p_n-1} Y_{l,s}(u) \qquad s = 1,2$$

with $\begin{cases} Y_{l,1}(u) &= \dfrac{1}{q_n}\sum_{i=2lq_n+1}^{(2l+1)q_n}[u(X_i^*) - \mathbb{E}(u(X_i^*))], \\ Y_{l,2}(u) &= \dfrac{1}{q_n}\sum_{i=2(2l+1)q_n+1}^{(2l+2)q_n}[u(X_i^*) - \mathbb{E}(u(X_i^*))]. \end{cases}$

To bound $P(\bar\nu_{n,1}(\varphi_\lambda\varphi_{\lambda'}) \geq B_{\lambda,\lambda'}x + V_{\lambda,\lambda'}\sqrt{2\|f\|_{\infty,A_1}x})$, we will use the Bernstein inequality given in Birgé and Massart (1998). A fast computation gives $\mathbb{E}|Y_{l,1}(\varphi_\lambda\varphi_{\lambda'})|^p \leq 2^{p-2}(B_{\lambda,\lambda'})^{p-2}(\sqrt{\|f\|_{\infty,A_1}}V_{\lambda,\lambda'})^2$. And then

$$P(|\bar\nu_{n,s}(\varphi_\lambda\varphi_{\lambda'})| \geq B_{\lambda,\lambda'}x + V_{\lambda,\lambda'}\sqrt{2\|f\|_{\infty,A_1}x}) \leq 2e^{-p_nx}.$$

Let $C = f_0^2[48\|f\|_{\infty,A_1}]^{-1}$, so that $x = 2C/L(\varphi)$. Given that $P(\Delta^c \cap \Omega_X^*) \leq P(D^c) \leq \sum_{\lambda,\lambda'} P\left(|\bar\nu_n(\varphi_\lambda\varphi_{\lambda'})| > B_{\lambda,\lambda'}x + V_{\lambda,\lambda'}\sqrt{2\|f\|_{\infty,A_1}x}\right)$,

$$P(\Delta^c \cap \Omega_X^*) \leq 4\mathcal{D}^2\exp\left\{-\frac{2p_nC}{L(\varphi)}\right\} \leq 4n^{1/(2\gamma+1)}\exp\left\{-C\frac{n}{q_nL(\varphi)}\right\}.$$

But $L(\varphi) \le \Phi_0 \mathcal{D}^2 \le \Phi_0 n^{1/(2\gamma+1)}$ and $q_n \le n^{1/2}$ so

$$P(\Delta^c \cap \Omega_X^*) \le 4n^{1/(2\gamma+1)} \exp\left\{-\frac{C}{\Phi_0} n^{\frac{2\gamma-1}{2(2\gamma+1)}}\right\} \le \frac{C'}{n^2}$$

because $\gamma > 1/2$.

### 6.6. Proof of Proposition 3

First we need to isolate even terms from odd terms in $Z_n^{(1)}(T)$ to avoid overlaps: $Z_n^{(1)}(T) = \frac{1}{2}(Z_n^{(1,1)}(T) + Z_n^{(1,2)}(T))$ with

$$\begin{cases} Z_n^{(1,1)}(T) = \dfrac{1}{n} \displaystyle\sum_{i=1, i \text{ odd}}^{n} V_T(Y_i, Y_{i+1}) - \mathbb{E}_X[V_T(Y_i, Y_{i+1})] \\ Z_n^{(1,2)}(T) = \dfrac{1}{n} \displaystyle\sum_{i=1, i \text{ even}}^{n} V_T(Y_i, Y_{i+1}) - \mathbb{E}_X[V_T(Y_i, Y_{i+1})] \end{cases}$$

It is sufficient to deal with the first term only, as the second one is similar. For each $i$, let $U_i = (Y_{2i-1}, Y_{2i})$, then

$$Z_n^{(1,1)}(T) = \frac{1}{n/2} \sum_{i=1}^{n/2} \{V_T(U_i) - \mathbb{E}_X[V_T(U_i)]\}.$$

Notice that conditionally to $X_1, \dots, X_n$, the $U_i$'s are independent. Thus we can use the Talagrand inequality recalled in Lemma 6 to bound

$$\mathbb{E}\left(\left[\sup_{T \in B_f(m,m')} Z_n^{(1,1)}(T)^2 - p_1(m, m')\right]_+\right).$$

We first remark that Property P1 entails $B_f(m, m') \subset \mathbb{S}_{m''}$ with $m'' = \max(m, m')$. Then, if $T$ belongs to $B_f(m, m')$,

$$T(x, y) = \sum_{j=J}^{m''} \sum_{kl} a_{jkl} \varphi_{jk}(x) \varphi_{jl}(y)$$

with $\sum_{jkl} a_{jkl}^2 = \|T\|^2 \le f_0^{-1}$.

• Let us bound $\|V_T\|_\infty$ for $T$ in $B_f(m, m')$. If $T(x, y) = \sum_{jkl} a_{jkl} \varphi_{jk}(x) \varphi_{jl}(y)$,

$$|V_T(x, y)|^2 \le \sum_{jkl} a_{jkl}^2 \sum_{jkl} |V_{\varphi_{jk} \otimes \varphi_{jl}}(x, y)|^2.$$

Then, since $V_{s \otimes t}(x, y) = v_s(x) v_t(y)$,

$$\sup_{T \in B_f(m,m')} |V_T(x, y)|^2 \le f_0^{-1} \sum_{jkl} |v_{\varphi_{jk}}(x) v_{\varphi_{jl}}(y)|^2.$$

But, according to Property P4, $\| \sum_k |v_{\varphi_{jk}}|^2 \|_\infty \leq \Phi_1 (2^j)^{2\gamma+2}$. So, using Lemma 4,

$$\sup_{T\in B_f(m,m')} \|V_T\|_\infty^2 \leq f_0^{-1}\Phi_1^2 \sum_{j=J}^{m''} (2^j)^{4\gamma+4} \leq f_0^{-1}\Phi_1^2 \frac{2^{4\gamma+4}}{2^{4\gamma+4}-1} D_{m''}^{4\gamma+4}$$

and $M_1 = f_0^{-1/2}\Phi_1 \sqrt{2^{4\gamma+4}/(2^{4\gamma+4}-1)} D_{m''}^{2\gamma+2}$.

- To compute $H^2$ we write

$$\begin{aligned}
&\mathbb{E}_X[\sup_{T\in B_f(m,m')} Z_n^{(1,1)}(T)^2] \leq f_0^{-1} \sum_{jkl} \mathbb{E}_X[Z_n^{(1,1)}(\varphi_{jk}\otimes\varphi_{jl})^2] \\
&\leq f_0^{-1} \sum_{jkl} \mathrm{Var}_X \left( \frac{1}{n} \sum_{i=1, i \text{ odd}}^{n} v_{\varphi_{jk}}(Y_i) v_{\varphi_{jl}}(Y_{i+1}) \right) \\
&\leq f_0^{-1} \sum_{jkl} \frac{1}{n} \mathrm{Var}_X \left( v_{\varphi_{jk}}(Y_1) v_{\varphi_{jl}}(Y_2) \right) \\
&\leq \frac{f_0^{-1}}{n} \sum_{jkl} \mathbb{E}_X[|v_{\varphi_{jk}}(Y_1)|^2 |v_{\varphi_{jl}}(Y_2)|^2] \qquad (16)
\end{aligned}$$

Here $\mathrm{Var}_X$ denotes the variance conditionally to $X_1, \ldots, X_{n+1}$. Now, for any function $G$, the following relation holds

$$\begin{aligned}
&\mathbb{E}_X[|G|^2(Y_1,Y_2)] = \mathbb{E}_X[|G|^2(X_1+\varepsilon_1, X_2+\varepsilon_2)] \\
&= \iint |G|^2(X_1+z_1, X_2+z_2) q(z_1) q(z_2) dz_1 dz_2 \\
&= \iint |G|^2(u_1,u_2) q(u_1-X_1) q(u_2-X_2) du_1 du_2 \leq \|q\|_\infty^2 \|G\|^2
\end{aligned}$$

Now, coming back to (16),

$$\mathbb{E}_X \left[ \sup_{T\in B_f(m,m')} Z_n^{(1,1)}(T)^2 \right] \leq \frac{f_0^{-1}}{n} \|q\|_\infty^2 \sum_{jkl} \|v_{\varphi_{jk}} \otimes v_{\varphi_{jl}}\|^2$$

$$\leq \frac{f_0^{-1}\|q\|_\infty^2}{n} \sum_j \left( \sum_k \|v_{\varphi_{jk}}\|^2 \right)^2 \leq \frac{\Phi_1^2 f_0^{-1} \|q\|_\infty^2}{n} \sum_{j=J}^{m''} (2^j)^{4\gamma+2},$$

using P5. Then, according to Lemma 4, $H^2 = \Phi_1^2 f_0^{-1}\|q\|_\infty^2 2^{4\gamma+2}/(2^{4\gamma+2} - 1)\dfrac{D_{m''}^{4\gamma+2}}{n}$.

- There remains to find $v$. First

$$\mathrm{Var}_X(V_T(Y_i,Y_{i+1})) \leq \mathbb{E}_X|V_T(Y_i,Y_{i+1})|^2 \leq \|q\|_\infty^2 \|V_T\|^2$$

We now observe that $\|V_T\|^2 = \|V_T^*\|^2/(4\pi^2)$ and then

$$\begin{aligned}
\|V_T\|^2 &= \frac{1}{4\pi^2}\iint \left|\frac{T^*(u,v)}{q^*(-u)q^*(-v)}\right|^2 dudv \\
&\leq \frac{1}{4\pi^2}\sqrt{\iint \frac{|T^*(u,v)|^2}{|q^*(-u)q^*(-v)|^4}dudv}\sqrt{\iint |T^*(u,v)|^2 dudv} \\
&\leq \frac{1}{4\pi^2}\sqrt{\sum_{jkl} a_{jkl}^2 \sum_{jkl}\iint \frac{|\varphi_{jk}^*(u)\varphi_{jl}^*(v)|^2}{|q^*(-u)q^*(-v)|^4}dudv}\sqrt{4\pi^2\|T\|^2}
\end{aligned}$$

For $T \in B_f(m,m')$,

$$\|V_T\|^2 \leq \frac{f_0^{-1/2}}{2\pi}\sqrt{f_0^{-1}\sum_j\sum_{kl}\int \frac{|\varphi_{jk}^*(u)|^2}{|q^*(-u)|^4}du\int \frac{|\varphi_{jl}^*(u)|^2}{|q^*(-u)|^4}du.}$$

But $(\varphi_{jk})^*(u) = 2^{-j/2}e^{iuk/2^j}\varphi^*(u/2^j)$ and then

$$\begin{aligned}
&\int \frac{|\varphi_{jk}^*(u)|^2}{|q^*(-u)|^4}du \leq \int \frac{2^{-j}|\varphi^*(u/2^j)|^2}{|q^*(-u)|^4}du \\
&\leq \int \frac{|\varphi^*(v)|^2}{|q^*(-v2^j)|^4}dv \leq k_0^{-4}(2^j)^{4\gamma}\int |\varphi^*(v)|^2(v^2+1)^{2\gamma}dv
\end{aligned}$$

Since $r > 2\gamma + 1/2$, Lemma 5 gives

$$\sum_{(k,l)\in\Lambda_j}\int \frac{|\varphi_{jk}^*(u)|^2}{|q^*(-u)|^4}du\int \frac{|\varphi_{jl}^*(u)|^2}{|q^*(-u)|^4}du \leq 3.2^{2j}C_{2,4\gamma}^2 k_0^{-8}(2^j)^{8\gamma}$$

Then, using Lemma 4 with $\rho = 8\gamma + 2$,

$$\|V_T\|^2 \leq \frac{C_{2,4\gamma}f_0^{-1}k_0^{-4}}{2\pi}\sqrt{\sum_{j=J}^{m''}3(2^j)^{8\gamma+2}} \leq \frac{C_{2,4\gamma}f_0^{-1}k_0^{-4}}{2\pi}\left(\frac{3.2^{8\gamma+2}}{2^{8\gamma+2}-1}\right)^{1/2} D_{m''}^{4\gamma+1}$$

and $v = \|q\|_\infty^2 C_{2,4\gamma}f_0^{-1}k_0^{-4}\sqrt{3.2^{8\gamma+2}}D_{m''}^{4\gamma+1}/(2\pi\sqrt{2^{8\gamma+2}-1})$.

We can now apply inequality (19)

$$\begin{aligned}
\mathbb{E}[\sup_{T\in B_f(m,m')}|Z_n^{(1,1)}(T)|^2 - 6H^2]_+ &\leq C\left(\frac{v}{n}e^{-k_1\frac{nH^2}{v}} + \frac{M_1^2}{n^2}e^{-k_2\frac{nH}{M_1}}\right) \\
&\leq C'\left(\frac{D_{m''}^{4\gamma+1}}{n}e^{-k_1' D_{m''}} + \frac{D_{m''}^{4\gamma+4}}{n^2}e^{-k_2'\sqrt{n}/D_{m''}}\right).
\end{aligned}$$

Yet there exists a positive constant $K$ such that

$$\sum_{m'\in\mathcal{M}_n} D_{m''}^{4\gamma+1}e^{-k_1' D_{m''}} \leq K.$$

Moreover, since $D_{m''} \leq n^{\frac{1}{4\gamma+2}}$, $D_{m''}^{4\gamma+4} e^{-k_2' \sqrt{n}/D_{m''}}/n \leq n^{1/(2\gamma+1)} e^{-k_2' n^{\gamma/(2\gamma+1)}}$ so that

$$\sum_{m' \in \mathcal{M}_n} D_{m''}^{4\gamma+4} e^{-k_2' \sqrt{n}/D_{m''}}/n^2 \leq K'/n.$$

Then, setting $K_1 = 6\Phi_1^2 f_0^{-1} \|q\|_\infty^2 2^{4\gamma+2}/(2^{4\gamma+2}-1)$,

$$\sum_{m' \in \mathcal{M}_n} \mathbb{E}[\sup_{T \in B_f(m,m')} |Z_n^{(1,1)}(T)|^2 - K_1 \frac{D_{m''}^{4\gamma+2}}{n}]_+ \leq \frac{C''}{n}$$

and the proposition is proved.

### *6.7. Proof of Proposition 4*

Since $\Pi_m$ belongs to $\mathbb{S}_m$, it can be written

$$\Pi_m(x,y) = \sum_{j'=J}^{m} \sum_{(k',l') \in \Lambda_{j'}} b_{j'k'l'} \varphi_{j'k'}(x) \varphi_{j'l'}(y)$$

with $\sum_{j'k'l'} b_{j'k'l'}^2 = \|\Pi_m\|^2 \leq \|\Pi\|_A^2$. From the embedding $B_f(m,m') \subset S_{m''}$ (where $m'' = \max(m,m')$), we have, if $T$ belongs to $B_f(m,m')$,

$$T(x,y) = \sum_{j=J}^{m''} \sum_{(k,l) \in \Lambda_j} a_{jkl} \varphi_{jk}(x) \varphi_{jl}(y)$$

with $\sum_{jkl} a_{jkl}^2 = \|T\|^2 \leq f_0^{-1}$.

We use the Talagrand inequality (19) in Lemma 6. But the variables $Y_i$ are not independent. We shall use the following approximation variables

$$\forall 1 \leq i \leq n+1 \qquad Y_i^* = X_i^* + \varepsilon_i.$$

These variables have the same properties as regards the $Y_i$'s as the $X_i^*$'s as regards the $X_i$'s (see (11)). More precisely, let, for $l = 0, \dots, p_n - 1$, $C_l = (Y_{2lq_n+1}, ..., Y_{(2l+1)q_n})$, $D_l = (Y_{(2l+1)q_n+1}, \dots, Y_{(2l+2)q_n})$, $C_l^* = (Y_{2lq_n+1}^*, ..., Y_{(2l+1)q_n}^*)$, $D_l^* = (Y_{(2l+1)q_n+1}^*, ..., Y_{(2l+2)q_n}^*)$. Then, since $A_l$ and $A_l^*$ have the same distribution and the sequences $(\varepsilon_i)$ and $(X_i)$ are independent, $C_l$ and $C_l^*$ have the same distribution. Moreover the construction of $A_l^*$ via Berbee's coupling Lemma implies that $C_l^*$ and $C_{l'}^*$ are independent if $l \neq l'$. At last $P(C_l \neq C_l^*) \leq \beta_{q_n}$.

Now we split $Z_n^{(2)}$ into two terms: $Z_n^{(2)}(T)\mathbb{1}_\Omega = (1/2)Z_n^{(2,1)}(T) + (1/2)Z_n^{(2,2)}(T)$ where

$$\begin{cases} Z_n^{(2,1)}(T) = \dfrac{1}{p_n} \displaystyle\sum_{l=0}^{p_n-1} \dfrac{1}{q_n} \sum_{i=2lq_n+1}^{(2l+1)q_n} Q_{T\Pi_m}(Y_i^*) - \mathbb{E}[Q_{T\Pi_m}(Y_i^*)] \\ Z_n^{(2,2)}(T) = \dfrac{1}{p_n} \displaystyle\sum_{l=0}^{p_n-1} \dfrac{1}{q_n} \sum_{i=(2l+1)q_n+1}^{(2l+2)q_n} Q_{T\Pi_m}(Y_i^*) - \mathbb{E}[Q_{T\Pi_m}(Y_i^*)] \end{cases}$$

Then we apply Talagrand's inequality to $Z_n^{(2,1)}(T)$.

• Let us first compute $M_1$. We have to bound $\|Q_{T\Pi_m}\|_\infty$ for $T$ in $B_f(m,m')$. By linearity of $Q$

$$Q_{T\Pi_m}(x)=\sum_{jkl}a_{jkl}\sum_{j'k'l'}b_{j'k'l'}Q_{\varphi_{jk}\varphi_{j'k'}\otimes\varphi_{jl}\varphi_{j'l'}}(x)$$

Then, since $Q_{s\otimes t}(x)=v_s(x)\int t(y)dy$, using the Schwarz inequality,

$$\begin{aligned}|Q_{T\Pi_m}(x)|^2 &\le \sum a_{jkl}^2 b_{j'k'l'}^2\sum|v_{\varphi_{jk}\varphi_{j'k'}}(x)\int\varphi_{jl}\varphi_{j'l'}|^2\\ &\le f_0^{-1}\|\Pi\|_A^2\sum_{jkk'l}|v_{\varphi_{jk}\varphi_{jk'}}(x)|^2\end{aligned}$$

since the $\varphi_{jl}$ are orthonormal. The property P6 gives then

$$\|Q_{T\Pi_m}\|_\infty^2\le f_0^{-1}\|\Pi\|_A^2\Phi_1\sum_{j=J}^{m''}(2^j)^{2\gamma+3}2^j$$

so that (using Lemma 4 again) $M_1=f_0^{-1/2}\|\Pi\|_A\sqrt{\Phi_1 2^{2\gamma+4}/(2^{2\gamma+4}-1)}D_{m''}^{\gamma+2}$.

• Now, we compute $H^2$. For $T\in B_f(m,m')$,

$$|Z_n^{(2,1)}(T)|^2\le\sum_{jkl}a_{jkl}^2\sum_{jkl}|Z_n^{(2,1)}(\varphi_{jk}\otimes\varphi_{jl})|^2$$

Thus

$$\begin{aligned}\mathbb{E}\left[\sup_{T\in B_f(m,m')}Z_n^{(2,1)}(T)^2\right]&\le f_0^{-1}\sum_{jkl}\mathbb{E}\left[Z_n^{(2,1)}(\varphi_{jk}\otimes\varphi_{jl})^2\right]\\ &\le f_0^{-1}\sum_{jkl}\mathrm{Var}\left[\frac{1}{p_n}\sum_{l=0}^{p_n-1}\frac{1}{q_n}\sum_{i=2lq_n+1}^{(2l+1)q_n}Q_{(\varphi_{jk}\otimes\varphi_{jl})\Pi_m}(Y_i^*)\right]\end{aligned}$$

The variables $(C_l^*)$ are independent and identically distributed so

$$\mathbb{E}\left[\sup_{T\in B_f(m,m')}Z_n^{(2,1)}(T)^2\right]\le f_0^{-1}\sum_{jkl}\frac{1}{p_n}\mathrm{Var}\left[\frac{1}{q_n}\sum_{i=1}^{q_n}Q_{(\varphi_{jk}\otimes\varphi_{jl})\Pi_m}(Y_i^*)\right]$$

However, on $\Omega$, $C_1$ and $C_1^*$ have the same distribution, so that

$$\mathrm{Var}\left[\frac{1}{q_n}\sum_{i=1}^{q_n}Q_{(\varphi_{jk}\otimes\varphi_{jl})\Pi_m}(Y_i^*)\right]=\mathrm{Var}\left[\frac{1}{q_n}\sum_{i=1}^{q_n}Q_{(\varphi_{jk}\otimes\varphi_{jl})\Pi_m}(Y_i)\right].$$

And, coming back to the definition of $Q_T$, for $i_1\neq i_2$,

$$\begin{aligned}&\mathrm{Cov}(Q_{(\varphi_{jk}\otimes\varphi_{jl})\Pi_m}(Y_{i_1}),Q_{(\varphi_{jk}\otimes\varphi_{jl})\Pi_m}(Y_{i_2}))\\ &=\frac{1}{4\pi^2}\iint\mathbb{E}(e^{iY_{i_1}u}e^{-iY_{i_2}v})\frac{[(\varphi_{jk}\otimes\varphi_{jl})\Pi_m]^*(u,0)}{q^*(-u)}\frac{[(\varphi_{jk}\otimes\varphi_{jl})\Pi_m]^*(-v,0)}{q^*(v)}dudv\\ &=\frac{1}{4\pi^2}\iint\mathbb{E}(e^{iX_{i_1}u}e^{-iX_{i_2}v})[(\varphi_{jk}\otimes\varphi_{jl})\Pi_m]^*(u,0)[(\varphi_{jk}\otimes\varphi_{jl})\Pi_m]^*(-v,0)dudv\end{aligned}$$

since $\mathbb{E}(e^{i\varepsilon_{i_1}u}e^{-i\varepsilon_{i_2}v}) = q^*(-u)q^*(v)$. Now using (10),

$$\begin{aligned}&\mathrm{cov}(Q_{(\varphi_{jk}\otimes\varphi_{jl})\Pi_m}(Y_{i_1}), Q_{(\varphi_{jk}\otimes\varphi_{jl})\Pi_m}(Y_{i_2}))\\ &= \mathrm{cov}(\int(\varphi_{jk}\otimes\varphi_{jl})\Pi_m(X_{i_1},y)dy, \int(\varphi_{jk}\otimes\varphi_{jl})\Pi_m(X_{i_2},y)dy)\end{aligned}$$

It implies that

$$\begin{aligned}\mathrm{Var}\left[\frac{1}{q_n}\sum_{i=1}^{q_n} Q_{(\varphi_{jk}\otimes\varphi_{jl})\Pi_m}(Y_i)\right] &\leq \frac{1}{q_n^2}\sum_{i=1}^{q_n}\mathrm{Var}[Q_{(\varphi_{jk}\otimes\varphi_{jl})\Pi_m}(Y_i)]\\ &+\mathrm{Var}\left[\frac{1}{q_n}\sum_{i=1}^{q_n}\int(\varphi_{jk}\otimes\varphi_{jl})\Pi_m(X_i,y)dy\right]\end{aligned}$$

And then

$$\begin{aligned}\mathbb{E}\left[\sup_{T\in B_f(m,m')} Z_n^{(2,1)}(T)^2\right] &\leq f_0^{-1}\sum_{jkl}\frac{1}{p_nq_n}\mathrm{Var}[Q_{(\varphi_{jk}\otimes\varphi_{jl})\Pi_m}(Y_1)]\\ &+ f_0^{-1}\sum_{jkl}\frac{1}{p_n}\mathrm{Var}\left[\frac{1}{q_n}\sum_{i=1}^{q_n}\int(\varphi_{jk}\otimes\varphi_{jl})\Pi_m(X_i,y)dy\right] \quad (17)\end{aligned}$$

For the second term in (17), we use Lemma 7 to write

$$\begin{aligned}&\sum_{jkl}\mathrm{Var}\left[\frac{1}{q_n}\sum_{i=1}^{q_n}\int(\varphi_{jk}\otimes\varphi_{jl})\Pi_m(X_i,y)dy\right]\\ &\leq \frac{4\sum_k\beta_k}{q_n}\|\sum_{jkl}|\int(\varphi_{jk}\otimes\varphi_{jl})\Pi_m(.,y)dy|^2\|_\infty\end{aligned}$$

For all real $x$

$$\int(\varphi_{jk}\otimes\varphi_{jl})\Pi_m(x,y)dy = \sum_{k'} b_{jk'l}\varphi_{jk}\varphi_{jk'}(x)$$

Therefore

$$\sum_{jkl}|\int(\varphi_{jk}\otimes\varphi_{jl})\Pi_m(x,y)dy|^2 \leq \|\Pi\|_A^2\sum_{jkk'l}|\varphi_{jk}\varphi_{jk'}(x)|^2 \leq \|\Pi\|_A^2\Phi_1^2D_{m''}^3$$

using property P3. Then

$$\sum_{jkl}\mathrm{Var}\left[\frac{1}{q_n}\sum_{i=1}^{q_n}\int(\varphi_{jk}\otimes\varphi_{jl})\Pi_m(X_i,y)dy\right] \leq \frac{4\sum_k\beta_k}{q_n}\|\Pi\|_A^2\Phi_1^2D_{m''}^3$$

Thus we have bound the second term in (17) by $2f_0^{-1}\sum_k\beta_k\|\Pi\|_A^2\Phi_1^2D_{m''}^3/n$.

For the first term in (17), we bound $\sum_{jkl}\mathbb{E}[|Q_{(\varphi_{jk}\otimes\varphi_{jl})\Pi_m}(Y_1)|^2]$:

$$\begin{aligned}\sum_{jkl}\mathbb{E}[|Q_{(\varphi_{jk}\otimes\varphi_{jl})\Pi_m}(Y_1)|^2] &\le \sum_{j'k'l'} b^2_{j'k'l'}\sum_{jkl}\sum_{j'k'l'}|\int\varphi_{jl}\varphi_{j'l'}|^2\mathbb{E}|v_{\varphi_{jk}\varphi_{j'k'}}(Y_1)|^2\\ &\le \|\Pi\|_A^2\sum_{jkk'}2^j\mathbb{E}|v_{\varphi_{jk}\varphi_{jk'}}(Y_1)|^2\end{aligned}$$

But $\mathbb{E}|v_{\varphi_{jk}\varphi_{jk'}}(Y_1)|^2=\int|v_{\varphi_{jk}\varphi_{jk'}}(x)|^2p(x)dx$ where $p$ is the density of $Y_1$. Since $p=q*f$, $|p(x)|\le\|q\|_\infty$ for all $x$. Then

$$\mathbb{E}|v_{\varphi_{jk}\varphi_{jk'}}(Y_1)|^2\le\|q\|_\infty\int|v_{\varphi_{jk}\varphi_{jk'}}(x)|^2dx$$

and

$$\sum_{jkl}\mathbb{E}[|Q_{(\varphi_{jk}\otimes\varphi_{jl})\Pi_m}(Y_1)|^2]\le\|\Pi\|_A^2\|q\|_\infty\sum_{j=J}^{m''}2^j\sum_{kk'}\int|v_{\varphi_{jk}\varphi_{jk'}}(x)|^2dx$$

$$\le\|\Pi\|_A^2\|q\|_\infty\Phi_1\sum_{j=J}^{m''}(2^j)^{2\gamma+3}\le\|\Pi\|_A^2\|q\|_\infty\Phi_1\frac{2^{2\gamma+3}}{2^{2\gamma+3}-1}D_{m''}^{2\gamma+3},$$

applying Property P7 and Lemma 4. We finally obtain

$$\begin{aligned}\mathbb{E}\left[\sup_{T\in B_f(m,m')}Z_n^{(2,1)}(T)^2\right] &\le 2f_0^{-1}\|\Pi\|_A^2\|q\|_\infty\Phi_1\frac{2^{2\gamma+3}}{2^{2\gamma+3}-1}\frac{D_{m''}^{2\gamma+3}}{n}\\ &\quad+2f_0^{-1}\sum_k\beta_k\|\Pi\|_A^2\Phi_1^2\frac{D_{m''}^3}{n}\end{aligned}$$

Since the order of $nH^2$ has to be larger than the one of $v$, we choose

$$H^2=2f_0^{-1}\|\Pi\|_A^2\Phi_1\max(\|q\|_\infty 2^{2\gamma+3}/(2^{2\gamma+3}-1),\Phi_1)\left[\frac{D_{m''}^{2\gamma+7/2}}{n}+(\sum_k\beta_k)\frac{D_{m''}^3}{n}\right].$$

- Lastly, using Lemma 7 again,

$$\begin{aligned}&\mathrm{Var}(\frac{1}{q_n}\sum_{i=2lq_n+1}^{(2l+1)q_n}Q^2_{T\Pi_m}(Y_i^*))=\mathrm{Var}(\frac{1}{q_n}\sum_{i=2lq_n+1}^{(2l+1)q_n}Q^2_{T\Pi_m}(Y_i))\\ &\le\frac{4}{q_n}\mathbb{E}[|Q_{T\Pi_m}|^2(Y_1)b(Y_1)]\le\frac{4}{q_n}\|Q_{T\Pi_m}\|_\infty(\mathbb{E}[|Q_{T\Pi_m}|^2(Y_1)])^{1/2}(\mathbb{E}[b^2(Y_1)])^{1/2}\\ &\le\frac{4\sqrt{2\sum_k(k+1)\beta_k}}{q_n}\|Q_{T\Pi_m}\|_\infty(\mathbb{E}[|Q_{T\Pi_m}|^2(Y_1)])^{1/2}\qquad(18)\end{aligned}$$

We have already proved that $\|Q_{T\Pi_m}\|_\infty\le f_0^{-1/2}\|\Pi\|_A\sqrt{\Phi_1 2^{2\gamma+4}/(2^{2\gamma+4}-1)}D_{m''}^{\gamma+2}$.

Now we need a sharp bound on $\mathbb{E}[|Q_{T\Pi_m}(Y_1)|^2]$. We have

$$\mathbb{E}[|Q_{T\Pi_m}(Y_1)|^2] \leq \|q\|_\infty \int |Q_{T\Pi_m}|^2 = \frac{\|q\|_\infty}{2\pi} \int \left| \frac{(T\Pi_m)^*(u,0)}{q^*(-u)} \right|^2 du$$

Then it follows from the Schwarz inequality that

$$\mathbb{E}[|Q_{T\Pi_m}(Y_1)|^2] \leq \frac{\|q\|_\infty}{2\pi} \sqrt{\int \frac{|(T\Pi_m)^*(u,0)|^2}{|q^*(-u)|^4} du} \sqrt{\int |(T\Pi_m)^*(u,0)|^2 du}$$

We will evaluate the two terms under the square roots. First observe that

$$(T\Pi_m)^*(u,0) = \sum_{jkl} \sum_{j'k'l'} a_{jkl} b_{j'k'l'} (\varphi_{jk}\varphi_{j'k'})^*(u)(\varphi_{jl}\varphi_{j'l'})^*(0)$$
$$= \sum_{jkk'l} a_{jkl} b_{jk'l} (\varphi_{jk}\varphi_{jk'})^*(u)$$

since $(\varphi_{jl}\varphi_{j'l'})^*(0) = \int \varphi_{jl}\varphi_{j'l'} = \mathbb{1}_{j=j',l=l'}$. Then

$$\int |(T\Pi_m)^*(u,0)|^2 du \leq \sum_{jkk'l} a_{jkl}^2 b_{jk'l}^2 \sum_{jkk'l} \int |(\varphi_{jk}\varphi_{jk'})^*(u)|^2 du$$
$$\leq 2\pi f_0^{-1} \|\Pi\|_A^2 \sum_{jkk'l} \int |(\varphi_{jk}\varphi_{jk'})(u)|^2 du$$
$$\leq 2\pi f_0^{-1} \|\Pi\|_A^2 \sum_j 2^j \| \sum_{k'} |\varphi_{jk'}|^2 \|_\infty \sum_{k=k'-2N+2}^{k'+2N-2} \int |\varphi_{jk'}|^2$$

by taking into account the superposition of the supports. Using now (8)

$$\int |(T\Pi_m)^*(u,0)|^2 du \quad \leq \quad 2\pi f_0^{-1} \|\Pi\|_A^2 \sum_j 2^j C'(\varphi) 2^j (4N-3)$$
$$\leq \quad 2\pi f_0^{-1} \|\Pi\|_A^2 \Phi_1 (4N-3) \frac{2}{3} D_{m''}^2$$

Now

$$\int \frac{|(T\Pi_m)^*(u,0)|^2}{|q^*(-u)|^4} du \leq \sum_{jkk'l} a_{jkl}^2 b_{jk'l}^2 \sum_{jkk'l} \int \frac{|(\varphi_{jk}\varphi_{jk'})^*(u)|^2}{|q^*(-u)|^4} du$$
$$\leq f_0^{-1} \|\Pi\|_A^2 \sum_{jkk'l} \int 2^j \frac{|(\varphi_{jk}\varphi_{jk'})^*(2^j v)|^2}{|q^*(-2^j v)|^4} dv$$

Hence, inequality (9) and Assumption H1 show that

$$\int \frac{|(T\Pi_m)^*(u,0)|^2}{|q^*(-u)|^4} du$$
$$\leq f_0^{-1} \|\Pi\|_A^2 \sum_{jkk'l} \int 2^j C_r^2 \left[ |v|^{2(1-r)} \mathbb{1}_{|v|>1} + \mathbb{1}_{|v|\leq 1} \right] k_0^{-4} ((2^j v)^2 + 1)^{2\gamma} dv$$
$$\leq f_0^{-1} \|\Pi\|_A^2 C_r^2 k_0^{-4} C \sum_{jkk'l} (2^j)^{4\gamma+1}$$

with $C = \int \left[|v|^{2(1-r)}\mathbb{1}_{|v|>1} + \mathbb{1}_{|v|\le 1}\right](v^2+1)^{2\gamma}dv < \infty$ as soon as $r > 2\gamma + 3/2$. Then

$$\begin{aligned}\int \frac{|(T\Pi_m)^*(u,0)|^2}{|q^*(-u)|^4}du &\le f_0^{-1}\|\Pi\|_A^2 C_r^2 k_0^{-4} C \sum_{j=J}^{m}\sum_{k'l}\sum_{k=k'-2N+2}^{k'+2N-2}(2^j)^{4\gamma+1}\\ &\le f_0^{-1}\|\Pi\|_A^2 C_r^2 k_0^{-4} C3(4N-3)\frac{2^{4\gamma+3}}{2^{4\gamma+3}-1}D_{m''}^{4\gamma+3}\end{aligned}$$

Finally

$$\mathbb{E}[|Q_{T\Pi_m}(Y_1)|^2] \le \frac{\|q\|_\infty}{\pi} f_0^{-1}\|\Pi\|_A^2(4N-3)C_r k_0^{-2}\sqrt{C\frac{2^{4\gamma+3}}{2^{4\gamma+3}-1}\pi\Phi_1}D_{m''}^{2\gamma+5/2}$$

Then (18) gives

$$v = \sqrt{2\sum_k (k+1)\beta_k\|q\|_\infty} f_0^{-1}\|\Pi\|_A^2 k_0^{-1} C(\gamma,r,N,\Phi_1)\frac{D_{m''}^{2\gamma+13/4}}{q_n}.$$

Then replacing $n$ by $p_n$ in inequality (19) gives

$$\begin{aligned}\mathbb{E}[\sup_{T\in B_f(m,m')}|Z_n^{(2,1)}(T)|^2 - 6H^2]_+ &\le C\left(\frac{v}{p_n}e^{-k_1\frac{p_nH^2}{v}} + \frac{M_1^2}{p_n^2}e^{-k_2\frac{p_nH}{M_1}}\right)\\ &\le C'\left(\frac{D_{m''}^{2\gamma+13/4}}{n}e^{-k_1'D_{m''}^{1/4}} + \frac{D_{m''}^{2\gamma+4}q_n^2}{n^2}e^{-k_2'\frac{\sqrt{n}}{q_nD_{m''}^{1/4}}}\right)\end{aligned}$$

where $C'$ and $k_1'$ depend on $r, N, \gamma, \Phi_1, f_0, \|\Pi\|_A, \|q\|_\infty, \sum_k(k+1)\beta_k$ and $\sum_k \beta_k$. But there exists a positive constant $K$ such that

$$\sum_{m'\in\mathcal{M}_n} D_{m''}^{2\gamma+13/4}e^{-k_1'D_{m''}^{1/4}} \le K.$$

Moreover $D_{m''}^{1/4} \le n^{1/8}$ and $q_n \le n^c$ with $c + 1/8 < 1/2$, which implies

$$\sum_{m'\in\mathcal{M}_n}\mathbb{E}\left[\sup_{T\in B_f(m,m')}|Z_n^{(2,1)}(T)|^2 - K_2\|\Pi\|_A^2\left(\frac{D_{m''}^{2\gamma+7/2}}{n} + (\sum_k\beta_k)\frac{D_{m''}^3}{n}\right)\right]_+ \le \frac{C''}{n}$$

with $K_2 = 12f_0^{-1}\Phi_1\max(\|q\|_\infty 2^{2\gamma+3}/(2^{2\gamma+3}-1), \Phi_1)$. Thus, if $p_2(m,m') = p_2^{(1)}(m,m') + p_2^{(2)}(m,m')$ with $p_2^{(1)}(m,m') = K_2\|\Pi\|_A^2 D_{m''}^{2\gamma+7/2}/n$ and $p_2^{(2)}(m,m') = K_2\|\Pi\|_A^2\ (\sum_k\beta_k)D_{m''}^3/n$, then

$$\sum_{m'\in\mathcal{M}_n}\mathbb{E}\left(\left[\sup_{T\in B_f(m,m')}Z_n^{(2)}(T)^2 - p_2(m,m')\right]_+\mathbb{1}_\Omega\right) \le \frac{C_2}{n}.$$

### *6.8. Technical Lemmas*

**Lemma 4.** *For all $m \geq J$*

$$\sum_{j=J}^{m}(2^j)^\rho \leq \frac{2^\rho}{2^\rho-1}D_m^\rho$$

*Proof of Lemma 4:* It is sufficient to write

$$\sum_{j=J}^{m}(2^j)^\rho = \frac{2^{\rho(m+1)}-2^{\rho J}}{2^\rho-1} \leq \frac{2^\rho}{2^\rho-1}2^{\rho m} \leq \frac{2^\rho}{2^\rho-1}D_m^\rho$$

**Lemma 5.** *If $|\varphi^*(x)| \leq k_3(x^2+1)^{-r/2}$ for all real $x$ then*
- *if $s$ and $\alpha$ are reals such that $sr > \alpha+1$*

$$\int|\varphi^*(x)|^s(x^2+1)^{\alpha/2}dx \leq C_{s,\alpha} < \infty$$

- *if $r > 1$*

$$\int|\varphi^*(y)\varphi^*(x-y)|dy \leq C_r(|x|^{1-r}\mathbb{1}_{|x|>1}+\mathbb{1}_{|x|\leq 1})$$

*Proof of Lemma 5:*
- For the first point, it is sufficient to observe that the function $(x^2+1)^{(-rs+\alpha)/2}$ is integrable if $-rs+\alpha > -1$.
- By changing the variable ($y = xu$), we get

$$\begin{aligned}&\int|\varphi^*(y)\varphi^*(x-y)|dy = \int|\varphi^*(xu)\varphi^*(x(1-u))|xdu\\ &\leq \int_{|u|>1/3 \text{ and } |1-u|>1/3} k_3|xu|^{-r}k_3|x(1-u)|^{-r}|x|du\\ &\quad+\int_{|u|\leq 1/3}k_3^2|x(1-u)|^{-r}|x|du+\int_{|1-u|\leq 1/3}k_3^2|xu|^{-r}|x|du\end{aligned}$$

$$\begin{aligned}\int|\varphi^*(y)\varphi^*(x-y)|dy &\leq k_3^2 3^r|x|^{1-2r}\int_{|u|>1/3}\frac{du}{|u|^r}+k_3^2|x|^{1-r}\frac{2}{3}\left|\frac{3}{2}\right|^r\\ &\quad+k_3^2|x|^{1-r}\frac{2}{3}\left|\frac{3}{2}\right|^r\\ &\leq k_3^2\left[\frac{2.3^{2r-1}}{r-1}|x|^{1-2r}+2^{2-r}3^{r-1}|x|^{1-r}\right]\end{aligned}$$

Thus, if $|x|>1$, $\int|\varphi^*(y)\varphi^*(x-y)|dy \leq C_r|x|^{1-r}$ and if $|x|\leq 1$, $\int|\varphi^*(y)\varphi^*(x-y)|dy \leq C_r$ with $C_r = k_3^2(2.3^{2r-1}/(r-1)+2^{2-r}3^{r-1})$.

**Lemma 6.** *Let $T_1, \dots, T_n$ be independent random variables and*

$$\nu_n(r) = (1/n)\sum_{i=1}^{n}[r(T_i) - \mathbb{E}(r(T_i)],$$

*for $r$ belonging to a countable class $\mathcal{R}$ of measurable functions. Then, for $\epsilon > 0$,*

$$\mathbb{E}[\sup_{r\in\mathcal{R}} |\nu_n(r)|^2 - 6H^2]_+ \leq C\left(\frac{v}{n}e^{-k_1\frac{nH^2}{v}} + \frac{M_1^2}{n^2}e^{-k_2\frac{nH}{M_1}}\right) \tag{19}$$

*with $k_1 = 1/6$, $k_2 = 1/(21\sqrt{2})$ and $C$ a universal constant and where*

$$\sup_{r\in\mathcal{R}} \|r\|_\infty \leq M_1, \quad \mathbb{E}\left(\sup_{r\in\mathcal{R}} |\nu_n(r)|\right) \leq H, \quad \sup_{r\in\mathcal{R}} \frac{1}{n}\sum_{i=1}^{n} \mathrm{Var}(r(T_i)) \leq v.$$

Usual density arguments allow using this result with non-countable class of functions $\mathcal{R}$.

*Proof of Lemma 6:* We apply the Talagrand concentration inequality given in Klein and Rio (2005) to the functions $s^i(x) = r(x) - \mathbb{E}(r(T_i))$ and we obtain

$$P(\sup_{r\in\mathcal{R}} |\nu_n(r)| \geq H + \lambda) \leq \exp\left(-\frac{n\lambda^2}{2(v + 4HM_1) + 6M_1\lambda}\right).$$

Then we modify this inequality following Birgé and Massart (1998) Corollary 2 p.354. It gives

$$P(\sup_{r\in\mathcal{R}} |\nu_n(r)| \geq (1+\eta)H + \lambda) \leq \exp\left(-\frac{n}{3}\min\left(\frac{\lambda^2}{2v}, \frac{\min(\eta,1)\lambda}{7M_1}\right)\right).$$

To find inequality (19) we use the formula $\mathbb{E}[X]_+ = \int_0^\infty P(X \geq t)dt$ with $X = \sup_{r\in\mathcal{R}} |\nu_n(r)|^2 - 6H^2$.

**Lemma 7.** *(Viennet (1997)) Let $(T_i)$ a strictly stationary process with $\beta$-mixing coefficients $\beta_k$. Then there exists a function $b$ such that*

$$\mathbb{E}[b(T_1)] \leq \sum_k \beta_k \quad \textit{and} \quad \mathbb{E}[b^2(T_1)] \leq 2\sum_k (k+1)\beta_k$$

*and for all function $\psi$ (such that $\mathbb{E}[\psi^2(T_1)] < \infty$) and for all $N$*

$$\mathrm{Var}(\sum_{i=1}^{N} \psi(T_i)) \leq 4N\mathbb{E}[|\psi|^2(T_1)b(T_1)].$$

In particular, for functions $(\psi_\lambda)$, $\sum_\lambda \mathrm{Var}(\sum_{i=1}^N \psi_\lambda(T_i)) \leq 4N(\sum_k \beta_k)\| \sum_\lambda |\psi_\lambda|^2\|_\infty$.

## Acknowledgements

I would like to thank Fabienne Comte for her careful readings of this work and Erwan Le Pennec for his help concerning the simulations.